\documentclass{amsart}
\usepackage{amssymb,epsfig,amsmath,latexsym,amsthm}
\theoremstyle{plain} \newtheorem{theorem}{Theorem}[section]
\newtheorem{lemma}[theorem]{Lemma} 
\newtheorem{proposition}[theorem]{Proposition}
\newtheorem{corollary}[theorem]{Corollary}

\newtheorem{conjecture}[theorem]{Conjecture}
\newtheorem{problem}[theorem]{Problem}
\newtheorem{ComplexityConjecture}[theorem]{Hyperbolic Complexity
Conjecture} \theoremstyle{definition}
\newtheorem{SmallestManifoldConjecture}[theorem]{Smallest Hyperbolic
Manifold Conjecture} \theoremstyle{definition}
\newtheorem{definition}[theorem]{Definition}
\newtheorem{remark}[theorem]{Remark}

\newtheorem{remarks}[theorem]{Remarks}

\newtheorem{notation}[theorem]{Notation}
\newtheorem{example}[theorem]{Example}

\DeclareMathOperator{\base}{base}

\DeclareMathOperator{\valence}{valence}

\DeclareMathOperator{\bridges}{bridges}

\DeclareMathOperator{\beams}{beams}

\newcommand{\cirM}{\overset\circ{M}}
\newcommand{\cirN}{\overset\circ{N}}

\newcommand{\cirV}{\overset\circ{V}}

\newcommand{\cirX}{\overset\circ{X}}

  \newcommand{\Sinfty}{S^2_{\infty}}

\newcommand{\BZ}{\mathbb Z}

\newcommand{\Vol}{\mathrm{Vol}}

\begin{document}

\title{Mom Technology and volumes of hyperbolic 3-manifolds}

\author{David Gabai}\footnote{Partially supported by NSF grants
DMS-0071852, DMS-0346270, and DMS-0504110.}\address{Department
  of Mathematics\\Princeton University\\Princeton, NJ 08544}

\author{Robert Meyerhoff}\address{Department of Mathematics, Boston
  College, Chestnut Hill, MA 02467}

\author{Peter Milley}\address{Department of Mathematics, University of
Calfornia Riverside, Riverside, CA 92521}

\thanks{Version 0.80, July 27, 2006}

\maketitle

\setcounter{section}{-1}

\section{Introduction}\label{S0}

This paper is the first in a series whose goal is to understand the
structure of low-volume complete orientable hyperbolic 3-manifolds.
Here we introduce \emph{Mom technology} and enumerate the hyperbolic
Mom-$n$ manifolds for $n\le 4$. Our long-term goal is to show that all
low-volume closed and cusped hyperbolic 3-manifolds are obtained by
filling a hyperbolic Mom-$n$ manifold, $n\le 4$ and to enumerate the
low-volume manifolds obtained by filling such a Mom-$n$.

William Thurston has long promoted the idea that volume is a good
measure of the complexity of a hyperbolic 3-manifold
(see, for example, \cite{Th1} page 6.48).
Among known low-volume manifolds, Jeff
Weeks (\cite{We}) and independently Sergei Matveev and
Anatoly Fomenko (\cite{MF}) have observed that there is a close connection
between the volume of closed hyperbolic 3-manifolds
and combinatorial complexity.  One goal of this project is
to explain this phenomenon, which is summarized by the following:

\begin{ComplexityConjecture} (Thurston,
Weeks, Matveev-Fomenko) The complete low-volume hyperbolic 3-manifolds
can be obtained by filling cusped hyperbolic 3-manifolds of small
topological complexity.\end{ComplexityConjecture}

\begin{remark}  Part of the challenge of this conjecture is to 
clarify the undefined adjectives \emph{low} and \emph{small}.  In the
late 1970's, Troels Jorgensen proved that for any positive constant
$C$ there is a finite collection of cusped hyperbolic 3-manifolds from
which all complete hyperbolic 3-manifolds of volume less than or equal
to $C$ can be obtained by Dehn filling.  Our long-term goal stated
above would constitute a concrete and satisfying realization of
Jorgensen's Theorem for ``low'' values of $C$.
\end{remark}

A special case of the Hyperbolic Complexity Conjecture is the
long-standing
\begin{SmallestManifoldConjecture}The Weeks Manifold $M_W$, obtained
  by $(5,1)$, $(5,2)$ surgery on the two components of the Whitehead
  Link, is the unique oriented hyperbolic $3$-manifold of minimum
  volume.
\end{SmallestManifoldConjecture}
Note that the volume of $M_W$ is $0.942\ldots$.

All manifolds in this paper will be orientable and all hyperbolic
structures are complete.  We call a compact manifold \emph{hyperbolic}
if its interior supports a complete hyperbolic structure of
finite volume.

\begin{definition}  A \emph{Mom-n structure} $(M,T,\Delta)$ consists of a
compact 3-manifold $M$ whose boundary is a union of tori, a preferred
boundary component $T$, and a handle decomposition $\Delta$ of the
following type.  Starting from $T\times I$, $n$ 1-handles and $n$
2-handles are attached to $T\times 1$ such that each 2-handle goes
over exactly three 1-handles, counted with multiplicity.  Furthermore,
each 1-handle encounters at least two 2-handles, counted with
multiplicity. We say that $M$ is a \emph{Mom-n} if it possesses a
Mom-$n$ structure \((M,T,\Delta)\).
\end{definition}

\begin{remarks}  On a Mom-$n$, the handle decomposition $\Delta$ deformation
retracts to an almost simple 2-complex which has $2n$ true vertices,
in the sense of Matveev \cite{Mv2}. Therefore Mom-$n$ manifolds are a
subset of those with Matveev complexity at most $2n$.\end{remarks}

Here is the fundamental idea at the foundation of our project.  Given
a complete finite-volume hyperbolic 3-manifold $N$, start with either
a slightly shrunken maximal horotorus neighborhood $V$ of a cusp or
slightly shrunken maximal tube V about a geodesic.  After expanding
$V$ in the normal direction, it eventually encounters itself, thereby
creating a 1-handle.  Subsequent expansions give rise to the creation
of 1, 2, and 3-handles.  In the presence of low volume we expect that
$V$ will rapidly encounter 1 and 2-handles and $\partial V$ together
with a subset of these handles (perhaps somewhat perturbed to allow
for the ``valence-3 2-handle condition'') will create a Mom-$n$
manifold $M$, for some $n\le 4$. Furthermore, the complement of $M$
will consist of cusp neighborhoods and tubular neighborhoods of
geodesics. In practice, the handle structure may arise in a somewhat
different manner; e.g., as a sub-complex of the dual triangulation of
the Ford Domain (see \cite{GMM3}).

The papers \cite{GM} and \cite{GMM} can be viewed as steps in this
direction when $V$ is a tubular neighborhood about a geodesic
$\gamma$.  Indeed, \cite{GM} gives a lower bound on $\Vol(N)$ in
terms of the tube radius of $\gamma$ and \cite{GMM} gives a lower
bound in terms of the first two ortholengths, or equivalently the
radii of the expanding $V$ as it encounters its first and second
1-handles.

\begin{definition}  If $i:M\to N$ is an embedding, then we say that the
embedding is \emph{elementary} if $i_*\pi_1(M)$ is abelian, and
\emph{non-elementary} otherwise.
\end{definition}

In \S 1 we give the basic definitions regarding Mom-$n$ manifolds
embedded in hyperbolic 3-manifolds and state for later use some
standard results about hyperbolic 3-manifolds and embedded tori in
such.  The end result of \S 2 - \S 4 is that if $n\le 4$, then given a
non-elementary Mom-$n$ in a hyperbolic manifold one can find a
non-elementary hyperbolic Mom-$k$, where $k\le n$. Weaker results are
given for general values of $n$.  In \S5 we enumerate the hyperbolic
Mom-$n$ manifolds for $n\le 4$; Theorem
\ref{list of mom-2's and mom-3's} and Conjecture \ref{number of mom-4's}
together imply the following:

\begin{theorem}There are $3$ hyperbolic Mom-$2$ manifolds, 21 hyperbolic
Mom-$3$ manifolds (including the 3 hyperbolic Mom-$2$'s, which are also
Mom-$3$'s).
\end{theorem}

\begin{conjecture}There are 138 hyperbolic Mom-$4$ manifolds (including
the hyperbolic Mom-$2$'s and Mom-$3$'s, which are also Mom-$4$'s).
\end{conjecture}

In \S 6 we show that any non-elementary embedding of a hyperbolic
Mom-$n$ manifold $M$, $n\le 4$, into a compact hyperbolic manifold $N$
gives rise to an \emph{internal Mom-$n$ structure} on $N$, i.e.\ every
component of $\partial M$ either splits off a cusp of $N$ or bounds a
solid torus in $N$.

In \S 7 we give examples of internal Mom-2
structures on cusped hyperbolic 3-manifolds, including in particular
a detailed
exposition of one of our key motivating examples.

In future papers we will use the Mom Technology to directly address
the Hyperbolic Complexity Conjecture, and the Smallest Hyperbolic
Manifold Conjecture. Indeed, in \cite{GMM3} we will identify all
$1$-cusped hyperbolic $3$-manifolds with volume less than $2.7$ by showing
that all such manifolds possess internal Mom-$n$ structures with
$n\le 3$. This result, in combination with work of Agol-Dunfield (see
\cite{AST}), gives a lower bound of $0.86$ for the volume of an orientable
hyperbolic 3-manifold (but see the note below).  The Agol-Dunfield
result is an improvement of an earlier result of Agol which provides a
tool for controlling the volume of a hyperbolic $3$-manifold in terms of
the volume of an appropriate cusped hyperbolic $3$-manifold from which
it is obtained via Dehn filling; the improved version utilizes recent
work of Perelman.

This leads to three very promising directions towards the Smallest
Hyperbolic Manifold Conjecture.  Either, improve the technology of
\cite{GMM3} to identify the 1-cusped manifolds of volume less than $2.852$,
and then apply Agol-Dunfield.  Or, extend the tube radius results of
\cite{GMT} from $\log(3)/2$ to $0.566$ (the Agol-Dunfield volume bound
involves the radius of a solid tube around a short geodesic in the
original closed hyperbolic 3-manifold).
Or, extend the method of \cite{GMM3} to the closed case, thereby providing
an essentially self-contained proof of the Smallest Hyperbolic
Manifold Conjecture.

Note that all three approaches require an analysis of volumes of
hyperbolic 3-manifolds obtained by Dehn filling a Mom-2 or Mom-3
manifold.  These Dehn filling spaces have been extensively studied by
J. Weeks and others, and it is highly likely that all low-volume
manifolds in these Dehn filling spaces have been identified.  However,
some work will need to be done to bring these studies up to a suitable
level of rigor.

The authors wish to thank Morwen Thistlethwaite for his assistance in
identifying the hyperbolic manifolds in \S 5.

\begin{section} {Basic Definitions and Lemmas}\end{section}

\vskip 12pt

\begin{definition}  Let $M$ be a compact connected
3-manifold $M$ with $B\subset \partial M$ a compact surface which may
be either disconnected or empty.  A \emph{handle structure} $\Delta $
on $(M,B)$ is the structure obtained by starting with $B\times I$
adding a finite union of 0-handles, then attaching finitely many 1 and
2-handles to $B\times 1$ and the 0-handles. We call $B\times I$
(resp. $B\times I \cup$ 0-handles) the \emph{base}
(resp. \emph{extended base}) and say that the handle structure is
based on $B$.  The \emph{valence} of a 1-handle is the number of
times, counted with multiplicity, the various 2-handles run over it
and the \emph{valence} of a 2-handle is the number of 1-handles,
counted with multiplicity, it runs over.

Following the terminology of Schubert \cite{Sch} and Matveev
\cite{Mv1} we call the 0-handles, 1-handles and 2-handles
\emph{balls}, \emph{beams} and \emph{plates} respectively.  We call
\emph{islands} (resp.  \emph{bridges}) the intersection of the
extended base with beams (resp. plates) and the closed complement of
the islands and bridges in $B\times 1 \cup \partial$(0-handles) are
the \emph{lakes}.  We say that $\Delta$ is \emph{full} if each lake is
a disc.  If $B=\emptyset$, then we say that $\Delta$ is a
\emph{classical handle structure}.

Let $M$ be a compact 3-manifold with $\partial M$ a union of tori and
let $T$ be a component of $\partial M$.  We say that $(M,T,\Delta)$ is
a \emph{weak Mom-n} if $\Delta$ is a handle structure based on $T$
without 0-handles or 3-handles and has an equal number of 1 and 2-handles, such
that each 1-handle is of valence $\ge 2$ and each 2-handle is of
valence 2 or 3.  Furthermore, there are exactly $n$ 2-handles of valence 3.
A weak Mom-$n$ with no valence-2 2-handles is a \emph{Mom-n}.  A weak
Mom-$n$ is \emph{strictly weak} if there exists a valence-2 2-handle.
\end{definition}

The following is a well-known existence result stated in our language.

\begin{proposition}  A compact 3-manifold $M$ has a weak Mom-structure if
and only if $\partial M$ is a union of at least two tori.
\end{proposition}

\begin{proof}
If $M$ has a weak Mom-$n$ structure, then by definition all of its
boundary components are tori and there is at least one such boundary
component.  Further, because there are no 3-handles in $\Delta$, there
must be another (torus) boundary component.

The converse is not much more difficult.  In fact, if $M$ has at least
two boundary components, then it is standard to create it by starting
with a thickened boundary component and adding 1 and 2-handles where
the 1 and 2-handles are of valence $\ge 2$.  By subdividing the
2-handles with 1-handles we satisfy the condition that the 2-handles
have valence $\le 3$.  Since $\chi(M)=0$, there are an equal number of
1 and 2-handles.\end{proof}

\begin{definition}  Call a torus that bounds a solid torus a
\emph{tube} and call a torus bounding a tube with knotted hole a \emph
{convolutube}.  Recall that a \emph{tube with knotted hole} is a
$B^3-\cirN(\gamma)$, where $\gamma$ is a knotted proper arc.
\end{definition}

The following standard result follows from the loop theorem (see for
example \cite{Jac}).

\begin{lemma}\label{torus}  If $S$ is a torus in an irreducible
3-manifold $N$, then either $S$ is incompressible or $S$ is a tube or
a convolutube.  If $S\subset \cirN$, $\partial N$ is incompressible
and $S$ is compressible, and there exists an embedded essential
annulus connecting $S$ to a component of $\partial N$, then $S$ is a
tube.  \qed\end{lemma}

\begin{proposition}  \label{convoluted embedding}If $M$ is a
non-elementary compact, connected 3-manifold embedded in the compact
hyperbolic 3-manifold $N$ and $\partial M$ is a union of tori, then,
up to isotopy, $N$ is obtained from $M$ by first filling a subset of
the components of $\partial M$ by solid tori to obtain the manifold
$M_1$, then compressing a subset of the components of $\partial M_1$
to obtain the manifold $M_2$, then attaching 3-balls to the 2-sphere
components of $\partial M_2$ to obtain $M_3$.  Furthermore all of
these operations can be performed within $N$.\end{proposition}

\begin{proof}  The components
of $\partial M$ that bound solid tori in $N$ are exactly those
boundary components which compress to the non-$M$ side.  Fill in all
such tori to obtain the manifold $M_1$.  If $P$ is a component of
$\partial M_1$ which is not boundary parallel in $N$, then $P$ is
compressible in $N$ and hence is a convolutube.  These convolutubes
can be isotoped to lie in pairwise disjoint 3-balls in $N$.  Therefore
we can compress all the compressible components of $\partial M_1$
(to obtain $M_2$) and
cap the resulting 2-spheres with 3-cells to obtain $M_3$ which is
isotopic to $N$.

Since $M_3$ must have all boundary components boundary parallel
in $N$ and $M_3$ is non-elementary, the result follows.\end{proof}

\begin{corollary}  Let $M\subset \cirN$ be a connected compact
non-elementary submanifold in the compact hyperbolic 3-manifold $N$.
If $\partial M$ is a union of tori, then each component of $\partial
N$ is parallel to a component of $\partial M$ via a parallelism
disjoint from $M$.\qed\end{corollary}

The following result is due to Kerckhoff (see \cite{Koj}).
\begin{lemma}\label{geodesic complement}  If $\gamma$ is a simple
closed geodesic in the complete, finite-volume hyperbolic 3-manifold
$N$, then $N-\gamma$ has a complete finite-volume hyperbolic
structure.\qed\end{lemma}

\begin{lemma}  \label{pi image}  Let $M$ be a compact submanifold of 
the compact hyperbolic 3-manifold $N$.

i) If $M_1=M-\cirX$, where $X\subset N$ is either a solid torus or 3-ball
with $\partial X\subset M$, or
$M_1$ is obtained by deleting a 2-handle or more generally
deleting an open regular neighborhood of a properly embedded arc from
M, then the inclusion $M\to N$ is a non-elementary embedding if and only
if the inclusion $M_1\to N$ is a non-elementary embedding.

ii) Suppose $M$ is non-elementary, $\partial M$ is a union of tori and $A$
is an essential annulus in $M$. Split $M$ along $A$; then some
component $M_1$ of the resulting manifold is non-elementary.

iii) Suppose $F\subset M$ is an embedded torus essential in $M$, and
$M\subset N$ is non-elementary. Split $M$ along $F$; then exactly one
component of the resulting manifold is non-elementary.
\end{lemma}

\begin{proof}  The conclusion is immediate in case (i) 
because both $M_1$ and $M$ have the same $\pi_1$-image.

Under the hypotheses of (ii) the boundary of each component of $M$
split along $A$ is also a union of tori.  We consider the case where
the split manifold connected, for the general case is similar.  Since all tori
in $N$ separate, $M$ is obtained from $M_1$ by attaching a thickened
annulus $A$ to a boundary parallel torus, a tube or a convolutube.  In
the first case $M $ and $M_1$ have the same $\pi_1$-image since $N$ is
anannular.  If $A$ is attached to the outside of a convolutube, then
$\pi_1(M)\subset \BZ$.  If $A$ is attached to the outside of a tube,
then $M$ can be enlarged to a Seifert fibered space in $N$ and hence
is elementary.

To prove (iii) note that $F$ is either boundary parallel or
is a tube or convolutube.  In each case the result follows
immediately.
\end{proof}

\begin{definition}  Let $N$ be a compact hyperbolic 3-manifold.
An \emph{internal Mom-n structure on N} consists of a non-elementary
embedding $f:M\to N$ where $(M,T,\Delta)$ is a Mom-$n$ and each
component of $\partial M$ is either boundary parallel in $N$ or bounds
a solid torus in $N$.  We will sometimes suppress mention of the
embedding and simply say that $(M,T,\Delta)$ is an internal Mom-$n$
structure on $N$.  In the natural way we define the notion of
\emph{weak internal Mom-$n$ structure on $N$}.
\end{definition}

\begin{lemma}  A non-elementary embedding of the Mom-$n$ manifold $M$ 
into the compact hyperbolic 3-manifold $N$ will fail to give an
internal Mom-$n$ structure on $N$ if and only if some component of
$\partial M$ maps to a convolutube.  In that case, a reimbedding of
$M$, supported in a neighborhood of the convolutubes gives rise to an
internal Mom-$n$ structure on $N$.\qed\end{lemma}

\begin{definition}  A \emph{general based Mom-n} $(M,B,\Delta)$ 
consists of a compact manifold $M$ with $\partial M$ a union of tori,
$B\subset \partial M$ a compact codimension-0 submanifold of
$\partial M$ that is
$\pi_1$-injective in $\partial M$, and $\Delta$ a handle structure
based on $B$ without 0-handles such that every 1-handle is of
valence-$\ge 2$, every 2-handle is of valence-3 and there are exactly
$n$ of each of them.  A \emph{weak general based Mom-n} is as above with
$\Delta$ perhaps having $k\ge 0$ extra valence-2 2-handles in addition
to the $n$ valence-3 2-handles and hence has $k+n$ 1-handles.

A \emph{general based internal Mom-n structure on N} consists of a
non-elementary
embedding $f:M\to N$, where $N$ is a compact hyperbolic 3-manifold and
$(M,B,\Delta)$ is a general based Mom-$n$ structure.  Along similar
lines we have the notion of \emph{weak general based internal Mom-n
structure on N}.
\end{definition}

\begin{notation}  If $X$ is a space, then let $|X|$ denote the number of
components of $X$, $\cirX$ denote the interior of $X$ and $N(X)$ denote
a regular neighborhood of $X$.  If $\sigma$ is a 2-handle, then let
$\delta \sigma$ denote the lateral boundary, i.e. the closure of that
part of $\partial \sigma$ which does not lie in lower index handles.
If $b$ is a bridge which lies in the 2-handle $\sigma$, then define
$\delta b=b\cap\delta\sigma$.  \end{notation}

\begin{section} {Handle structures and Normal surfaces}\end{section}

We slightly modify Haken's \cite{Ha} theory of surfaces in
handlebodies to our setting.  We closely parallel the excellent
exposition given by Matveev in \cite{Mv1}.

\smallskip
\begin{definition}  Let $\Delta$ be a handle structure on $M$ based 
on $B\subset \partial M$.  A compact surface $F\subset M$ is called
\emph{normal} if
\begin{enumerate}
\item $F$ intersects each plate $D^2\times I$ in parallel copies of
the form $D^2\times$ pt.
\item Each component of the intersection of $F$ with a beam $D^2\times
I$ is of the form $\alpha\times I$, where $\alpha$ is a proper arc
whose endpoints are disjoint from $\delta$(bridges).  Furthermore,
each component of $D^2\times 0-\alpha$ intersects $\delta$(bridges) in
at least two points.
\item Each component $U$ of $F\cap B\times [0,1]\cup \text{0-handles}$
is $\pi_1$-injective in $B\times [0,1]\cup \text{0-handles}$.  If
$U\cap B\times 0\neq\emptyset$, then $U$ is a product disc or annulus,
i.e. The inclusion $(U,U\cap(B\times 0),U\cap(B\times
1),U\cap(\partial B\times I)) \to (B\times I, B\times 0, B\times 1,
\partial B\times I)$, can be relatively isotoped to a vertical
embedding.
\item If $U$ is a component of $F\cap B\times I$ with $F\cap B\times
0\neq\emptyset$, then $U$ is an essential vertical disc or annulus,
i.e. $U=\alpha\times I$ where $\alpha\times 0$ is either an essential
simple closed curve or essential proper arc in $B\times 0$.

\end{enumerate}\end{definition}

\begin{remark}  i)  For $F$ closed, the second condition can be 
restated by requiring that $\alpha$ intersect distinct components of
$D^2\times 0\cap (\bridges)$.  When $\partial F\neq\emptyset$, the
second condition implies that $F$ is locally efficient in that it
neither can be locally boundary compressed nor can its weight be
reduced via an isotopy supported in the union of a 2-handle and its
neighborhing 1-handles.

ii) Note that $\partial F$ lies in the union of the beams, lakes and
$B\times 0$.\end{remark}

\begin{lemma}(Haken) If $F$ is a compact, incompressible,
boundary-incompressible surface in a compact irreducible 3-manifold,
then $F$ is isotopic to a normal surface.\qed\end{lemma}

\begin{definition}\label{def:valence}
Let $\Delta$ be a handle structure on $M$ based 
on $B\subset \partial M$.
The \emph{valence} $v(b)$ of a beam (resp. 
plate) is the number of plates (resp. beams) that attach to it,
counted with multiplicity.  Define the \emph{complexity} $C(\Delta)$
to be $(\rho_1(\Delta),|\Delta^1|)$ lexicographically ordered, where
$\rho_1(\Delta)=\Sigma_{\beams b} \max(v(b)-2,0)$ and $|\Delta^1|$ is
the number of 1-handles.  In particular we have\\

\noindent$\rho_1$-\textbf{formula}:
$\rho_1(\Delta)=\sum_{\text{2-handles}\ \sigma}\valence(\sigma) -
2|\Delta^1|+ |\text{valence-1}\ \text{1-handles}| +
2|\text{valence-0}\ \text{1-handles}|$
\end{definition}


\begin{lemma}(\text{Matveev}) \label{reduction} Let $\Delta$ be a 
handle structure on $M$ based on $B$, $F\subset M$ a closed normal
surface and let $M'$ be $M$ split along $F$.  If each component of
$M'\cap B\times [0,1]$ disjoint from $B\times 0$ is a 3-ball, then
$M'$ has a handle structure $\Delta'$ based on $B$ with
$\rho_1(\Delta')= \rho_1(\Delta)$.
\end{lemma}

\begin{proof}  This follows almost exactly as in \S 3 and
\S 4 of \cite{Mv1}: $M'$ naturally inherits a
handle structure $\Delta_1$ from $\Delta$ as follows.  The surface $F
$ splits $B\times I$ into various submanifolds one of which is
homeomorphic to $B\times [0,1]$ with $B\times 0=B$.  All of the other
submanifolds which lie in $M'$ are 3-balls.  This new $B\times [0,1]$
becomes the base and the 3-balls become 0-handles.  The various 1 and
2-handles are split by F into 1 and 2-handles and as in \cite{Mv1},
$\rho_1(\Delta_1)=\rho_1(\Delta)$.\end{proof}

\begin{lemma} \label{no trivial 1-handles}
Given the handle structure $\Delta$ on $(M,B)$, if some 1-handle is
valence-1, then there exists another structure $\Delta_1$ on $(M,B)$
with $C(\Delta_1)<C(\Delta)$.\qed\end{lemma}

\begin{lemma}\label{Mom-n complexity}  If $(M,T,\Delta)$ is a Mom-$n$,
then $C(\Delta)=(n,n)$.\qed\end{lemma}

\begin{lemma}  \label{classical reduction}
Let $\Delta$ be a handle structure on $(M,B)$, $F\subset M$ a
connected separating normal surface and $M_1$ be the component of
$M-\cirN(F)$ which does not contain $B$.  If each component of $F\cap
B\times [0,1]$ is a disc, then $M_1$ has a classical handle structure
$\Delta_1$ with $\rho_1(\Delta_1)\le \rho_1(\Delta)$.  \end{lemma}

\begin{proof}  This follows as in the proof of
Lemma \ref{reduction} after noting that each component of $M_1\cap
B\times [0,1]$ is a 3-ball and these 3-balls correspond to the
0-handles of the induced handle structure on $M_1$.\end{proof}

\begin{lemma}  If $\partial M$ is a union of tori,
and $\Delta$ is a handle structure on $(M,T)$ with $T$ a component of
$\partial M$, then there exists a weak Mom-$n$ $(M,T,\Delta_1)$ with
$n\le \rho_1(\Delta_1)$.  \end{lemma}

\begin{proof}  First apply Lemma \ref{no trivial 1-handles}, then add
1-handles to subdivide the valence-$k$, $k\ge 4$, 2-handles into
valence-3 2-handles.  \end{proof}

\begin{definition}  If $B\neq\emptyset$ is a compact submanifold of
$\partial M$, then define $\mathrm{rank}_{\rho_1}(M,B)$ to be the least
$n$ such that there exists a handle decomposition $\Delta$ on $(M,B)$
with $\rho_1(\Delta)=n$.\end{definition}

\begin{problem}\label{problem} Is there an example of a compact 
hyperbolic 3-manifold $N$ with $T$ a component of $\partial N$ and $A$
an essential annulus in $T$ such that $\mathrm{rank}_{\rho_1}(N,A)<$
$\mathrm{rank}_{\rho_1}(N,T)$?
\end{problem}

\begin{section} {Estimates for the reduction of $\rho_1$ under splitting}

The main result of these next two sections is Theorem \ref{hyperbolic
Mom-n} which shows that if a hyperbolic 3-manifold $N$ has an internal
Mom-$n$ structure $(M,T,\Delta)$ with $\Delta$ full and $n\le 4$, then
it has an internal Mom-$k$ structure $(M_1,T_1,\Delta_1)$ where
$k\le n$, $\Delta_1$ is full, and $M_1$ is hyperbolic. 
If $n>4$, we obtain the similar
conclusion except that ``full'' is replaced by ``general based'' and hence
$T_1$ can be a union of tori and annuli.

As far as we know, transforming a structure based on an annulus lying
in $T_1$ to one based on the whole torus $T_1$ may require an increase in
$\rho_1$.  This issue is responsible for many of the technicalities of
this section and the next.  See Problem \ref{problem}.

\begin{lemma}\label{free edge}  Let $f:M\to N$ be a non-elementary
embedding of a compact connected manifold into a compact irreducible
3-manifold.  Suppose $\partial M$ is a union of tori, $T$ is a
component of $\partial M$ and that $\Delta$ is a handlebody structure
on $(M,T)$ without 0-handles such that each 2-handle is of valence
$\ge 3$.  Then there exists a non-elementary embedding $g:M'\to N$
with $\partial M'$ a union of tori and a handle structure $\Delta'$ on
$(M',T')$ with $T'=T$ such that
$$\rho_1(\Delta')+2|\text{valence-1 1-handles of
$\Delta$}|+3|\text{valence-0 1-handles of $\Delta$}|\le
\rho_1(\Delta).$$

If instead each 2-handle of $\Delta$ is of valence $\ge 2$ then we
have
$$\rho_1(\Delta')+|\text{valence-1 1-handles of
$\Delta$}|+2|\text{valence-0 1-handles of $\Delta$}|\le
\rho_1(\Delta).$$
\end{lemma}

\begin{proof}  Both assertions follow similarly by
induction on the number of 1-handles of $\Delta$.  If $\eta$ is a
valence-1 1-handle, then cancelling $\eta$ with its corresponding
2-handle creates a handle structure $\Delta_1$. The Lemma follows by
applying the $\rho_1$-formula to $\Delta_1$ and induction.  If $\eta$
is a valence-0 1-handle, then the manifold $M_1$ obtained by deleting
$\eta$ is connected and $\chi(M_1)=1$, hence has a 2-sphere boundary
component $S$ which bounds a 3-ball disjoint from $\cirM_1$.  Let
$(M_2,T,\Delta_2)$ be obtained from $(M,T,\Delta)$ by deleting $\eta$
as well as a 2-handle which faces $S$. $M_2$ is a non-elementary
embedding by Lemma \ref{pi image}.  Now apply the $\rho_1$-formula and
induction.\end{proof}

\begin{lemma}  \label{boundary reduction} Let $f:M\to N$ be a 
non-elementary embedding of a manifold into a compact hyperbolic
3-manifold $N$, where $\partial M$ is a union of tori.  Suppose that
$M$ has a full handle structure $\Delta$ without 0-handles based on
the component $T$ of $\partial M$ such that every 2-handle is of
valence $\ge 3$.  If either of the following are true then there
exists a non-elementary embedding $f:M' \to N$ with handle structure
$\Delta'$ on $(M',T')$ such that $\rho_1(\Delta')+2\le
\rho_1(\Delta)$, $T'=T$ and $\partial M'$ is a union of tori:

i) There exists a valence-1 2-handle $\sigma\subset N-\cirM$ that can
be added to $\Delta$.

ii) There exists a disc $D\subset\partial M$ such that $\partial D$ is
the union of two arcs $\alpha\cup\beta$, where $\beta$ lies in a lake
and $\alpha$ lies in a 2-handle $\lambda$.  Furthermore, within
$\lambda\cap \partial M$, $\alpha$ separates components of
$\lambda\cap (\mathrm{1-handles})$.\end{lemma}

\begin{proof}  By Lemma
\ref{free edge} we can assume that every 1-handle of $\Delta$ is of
valence $\ge 2$.  To prove (i) let $\Delta_1$ be the handle structure
on the manifold $M_1$ obtained by attaching $\sigma$ to $\Delta$ along
$\partial M$.  Let $\eta$ denote the 1-handle which $\sigma$ meets.
Let $\Delta_2$ and $M_2$ be obtained by deleting a 2-handle $\lambda
\neq \sigma$, which faces the resulting 2-sphere boundary component.
Let $\Delta_3$ be obtained by cancelling $\sigma$ and $\eta$.
Finally, in the usual way, reduce to a non-elementary $M_4$ with
structure $\Delta_4$ on $(M_4,T)$ whose 1 and 2-handles are of valence
$\ge 2$.  Applying the $\rho_1$-formula shows that
$\rho_1(\Delta_4)+2\le \rho_1(\Delta)$ unless, measured in $\Delta$,
$\valence(\eta)=\valence(\lambda)=3$ and $\lambda$ attaches to $\eta$
at least twice.  If $\lambda$ attaches to $\eta$ twice, then in the
passage from $\Delta_1$ to $\Delta_2$ delete a 2-handle $\lambda_1\neq
\lambda$ which faces the 2-sphere.  If $\lambda$ attaches to $\eta$
thrice, then either $\Delta$ is not full or $\eta$ is the unique
1-handle of $\Delta$.  In the latter case $M_2=T\times I$ which is
elementary, a contradiction.

ii) Under these hypotheses we can attach a 2-handle $\sigma\subset
N-\cirM$ to $\Delta$ such that either valence($\sigma$)=1 or
$\valence(\sigma)\le \valence(\lambda)+2$ and $\lambda$ faces the
resulting 2-sphere boundary component.  If valence$(\sigma)=1$, then
apply (i).  Otherwise let $\Delta_1$ be obtained by deleting $\lambda$
and apply the $\rho_1$-formula and if necessary Lemma \ref{free edge}.
\end{proof}

\begin{lemma}\label {0-handle reduction}  Let $M$ be a non-elementary
embedding of a compact 3-manifold into the compact hyperbolic
3-manifold $N$ with $\partial M$ a union of tori, and let $\Delta$ be a
handle structure of $M$ based on $R\subset \partial M$.  If there
exists a valence $\ge 3$ 1-handle $\eta$ of $\Delta$ which attaches to
a 0-handle $\zeta$, then there exists a non-elementary embedding
$M'\to N$, and a handle structure $\Delta'$ based on $R'\subset
\partial M'$ such that $\rho_1(\Delta')<\rho_1(\Delta)$.  Here either
$(M',R')=(M,R)$ or $M'=M-\cirV$ and $R'=R\cup\partial V$ where $V$ is
an embedded solid torus in $\cirM$.
\end{lemma}

\begin{proof}  If $\eta$ also attaches to either the base or a 0-handle
distinct from $\zeta$, then cancelling $\eta$ with $\zeta$ gives rise
to $\Delta'$ on $(M,R)$ with $\rho_1(\Delta')<\rho_1(\Delta)$.  If
$\eta$ attaches only to $\zeta$, then let $M'$ be obtained by deleting
$\cirV$, the open solid torus gotten by hollowing out $\zeta$ and
$\eta$.  Let $R'=R\cup \partial V$ and let $\Delta'$ be the induced
structure on $(M',R')$.\end{proof}

\begin{lemma} \label{reducible reduction} Let $N$ be a compact hyperbolic
3-manifold.  If $(M,T,\Delta)$ is a full internal weak Mom-n structure
on $N$ and $M$ is reducible then there exists $(M_p,T_p,\Delta_p)$, a
weak internal Mom-k structure on N such that
$k+2=\rho_1(\Delta_p)+2\le\rho_1(\Delta)=n$. If $T$ does not lie in a
3-cell, then $T_p=T$.\end{lemma}

\begin{proof}  We first consider the case that no reducing sphere in 
$M$ bounds a ball in $N$ containing $T$.  Let $F$ be a least-weight
normal reducing 2-sphere and note that $F\cap T\times I$ is a union of
discs.  Let $M_0$ and $M_0'$ be the components of $M$ split along
$F$. By Lemma \ref{pi image} and the irreducibility of $N$ exactly one
of $M_0$, $M_0'$ is non-elementary.  We let $M_0$ (resp. $M_0'$) denote
the non-elementary (resp. elementary) component with $\Delta_0$
(resp. $\Delta_0'$) its induced structure.  By hypothesis $T\subset
M_0$.

We show that $\rho_1(\Delta_0^\prime)>0$.  Let $X$ denote the union of
the islands and bridges of $\Delta$ and $Y'=X\cap M_0^\prime$.  If
$\rho_1(\Delta_0^\prime)=0$, then each component $A$ of $Y'$ is an
annulus.  If some component $A$ of $Y'$ is disjoint from the lakes of
$\Delta$, then $Y'$ would be disjoint from the lakes and so $F$ would
2-fold cover a projective plane, which contradicts the fact that $N$
is irreducible.  Therefore each component of $Y'$ has one boundary
component in a lake and one component in $\cirX$.  Since $\Delta$ is
full, this implies that $F$ is a boundary parallel 2-sphere and
contradicts the fact that $\partial M$ is a union of tori.

If $Y=X\cap M_0$ then a similar argument shows that some component of
$Y$ is not an annulus and furthermore some 2-handle $\sigma$ of
$\Delta_0$ faces $F$ and attaches to a valence-$\ge 3$ 1-handle.
Delete $\sigma$ from $\Delta_0$ to obtain $(M_1,T,\Delta_1)$ with
$\rho_1(\Delta_1)+1\le\rho_1(\Delta_0)$.  The standard simplifying
moves as in Lemma \ref{free edge} transform $(M_1,T,\Delta_1)$ to
$(M_2,T,\Delta_2)$, a weak Mom-k with
$k=\rho_1(\Delta_2)\le\rho_1(\Delta_1)$. Since
$\rho_1(\Delta_0)+\rho_1(\Delta_0^\prime)=\rho_1(\Delta)=n$ we have
$k+2=\rho_1(\Delta_p)+2\le\rho_1(\Delta)=n$.  If $M_2$ is reducible,
then split along an essential least-weight sphere $F_2$ which is normal with
respect to $\Delta_2$.  Retain the component which contains $T$ and do
the usual operations to obtain the weak Mom-$k_2$ $(M_3,T,\Delta_3)$
with $k_2\le k$ and $M_3$ non-elementary.  By Haken finiteness, this
procedure terminates in a finite number of steps, completing the proof
in this case.

\smallskip We now consider the case that $T$ is compressible in $M$.
Let $F$ be a least-weight compressing disc for $T$.  Note that $F\cap
T\times I$ consists of discs and a single annulus.  If $M_1$ is $M$
split along $F$ with induced handle structure $\Delta_1$ and $A$ is
$T$ split along $F\cap T$, then $\Delta_1$ is based on $A$.  Since $M$
is obtained by attaching a 1-handle to a spherical component of
$\partial M_1$, it follows that $M_1$ is non-elementarily embedded in
$N$.  Let $M_2$ be obtained by attaching 2-handles to $\partial M_1$
along the components of $\partial A$, then capping off the resulting
2-sphere which faces $A$ with a 3-cell.  The resulting manifold $M_2$
has two 2-sphere boundary components and the induced handle structure
$\Delta_2$ is classical.  Since $\Delta$ is full some 2-handle of
$\Delta_2$ facing a 2-sphere of $\partial M_2$ attaches to a 1-handle
of valence $\ge 3$.  Delete this 2-handle to obtain the non-elementary
$M_3$ with handle structure $\Delta_3$ which satisfies
$\rho_1(\Delta_3)+1\le \rho_1(\Delta)$.  Delete another 2-handle to
create $M_4$ and $\Delta_4$ such that $\partial M_4$ is a union of
tori.  Now apply Lemma \ref{0-handle reduction} to $\Delta_4$ to
create $(M_5,T_5,\Delta_5)$ so that $\rho_1(\Delta_5)+1\le
\rho_1(\Delta_4)$.  Finally cancel the extraneous 0-handles and
valence-1 1-handles to create the desired weak internal Mom-k
structure $(M_6,T_6,\Delta_6)$.

\smallskip
Let $F$ be a least-weight essential normal 2-sphere for $\Delta$.
Since we can assume that $T$ is incompressible in $M$, $F\cap T\times
I$ is a union of discs.  If $T$ does not lie in the ball bounded by
$F\subset N$, then proceed as in the first part of the proof.
Otherwise let $M_0$ and $M_0'$ be the components of $M$ split along
$F$, with $M_0'$ the component containing $T$. Since $\Delta$ is full
and some component of $T\times 1\cap M_0'$ is nonplanar, it follows
that $\rho_1(\Delta_0')\ge 1$ and hence
$\rho_1(\Delta_0)+1\le\rho_1(\Delta)$.  To complete the proof apply
Lemma \ref{0-handle reduction} to $\Delta_0$, delete a 2-handle to
create a manifold with torus boundary components and cancel low
valence handles as in the previous paragraph.
\end{proof}

\begin{remark} If $(M, T, \Delta)$ is a weak internal Mom-n structure 
on the compact hyperbolic manifold $N$, and $\partial M$ is
compressible in $M$, then $M$ is reducible.  Therefore, if $M$ is
irreducible and $\Delta$ has an annular lake $A$ that is homotopically
inessential in $T\times 1$, then the core of $A$ bounds a disc in
$T\times I$ which separates off a 3-cell in $M$.  Absorbing this
3-cell into $T\times I$ simplifies $\Delta$ and transforms $A$ into a
disc lake.

From now on we will assume that if a homotopically inessential lake
appears it is immediately removed via the above operation.
\end{remark}

\begin{definition}  If $\Delta$ is a handle structure on $M$, then 
the \emph{sheets} of $\Delta$ are the connected components of the
space $\mathcal{S}$ which is the union of the 2-handles and the
valence-2 1-handles.  So sheets are naturally thickened surfaces which
are attached to a 3-manifold along their thickened boundaries.  The
valence of a sheet is the number of times the boundary runs over
1-handles counted with multiplicity.
\end{definition}

\begin{lemma}  \label{sheet reduction} Let $N$ be a compact hyperbolic
3-manifold and $f:M\to N$ a non-elementary embedding where $\partial M$
is a union of tori.  Let $\Delta$ be a handle structure on $M$ with no
0-handles based on a component $T$ of $\partial M$ such that the
valence of each 2-handle is at least 3 and $T$ does not lie in a
3-cell of $N$.  If some sheet $S$ of $\Delta$ is not a thickened disc
then there exists a non-elementary $M'\subset N$ with handle structure
$\Delta'$ based on a component $T'$ of $\partial M'$ such that
$\rho_1(\Delta')+1\le \rho_1(\Delta)$.  If equality must hold then
$S$ is a thickened annulus or Mobius band and if $\Delta$ is full
then $S$ is a thickened Mobius band and $\Delta'$ is full.
\end{lemma}

\begin{proof}  By Lemma \ref{free edge} and the proof of Lemma 
\ref{reducible reduction} we can assume that $M $ is irreducible and
  each 1-handle of $\Delta$ has valence $\ge 2$. If $M_1$ denotes the
  manifold obtained by deleting the sheet $S$, then
  $\chi(M)=\chi(M_1)+\chi(S)$.  Since $\chi(M)=0$, if $\chi(S)<0$,
  then $\partial M_1$ contains a 2-sphere.  This implies that $M$ is
  reducible.  Note that if $\chi(S)=0$, then $M_1$ is non-elementary.

Now assume that $\chi(S)=0$.  In this case $S$ is either an
annulus$\times I$ or a non-trivial $I$-bundle over a Mobius band.
Since $\chi(S)=0$, if $S$ contains more than one 2-handle, then
$\valence(S)\ge 2$.  If $\valence(S)>1$, then $\rho_1(\Delta_1)+2\le
\rho_1(\Delta)$ holds where $\Delta_1$ is the induced structure on
$M_1$.  If $S$ is a valence-1 annulus$\times I$, then $\Delta$ has an
annular lake and $\rho_1(\Delta_1)+1=\rho_1(\Delta)$.

If $S$ is a thickened Mobius band of valence-1, then
$\rho_1(\Delta_1)+1= \rho_1(\Delta)$.  If $S$ attached to the
component $R$ of $\partial M_1$ and $R$ had an annular lake, then $R$
must be a tube with a compressing disc $D$ whose boundary goes over a
1-handle $\eta$ of $\Delta_1$ exactly once. Let $M_2$ be obtained by
attaching the 2-handle $\sigma$ with core $D$ to $M_1$ and $\Delta_2$
the induced handle structure.  Proceed as in the proof of Lemma \ref
{boundary reduction} to show that exists a non-elementary embedding
$(M_3,T)$ with handle structure $\Delta_3$ with
$\rho_1(\Delta_3)+2\le\rho_1(\Delta_1)+1\le \rho_1(\Delta)$.
\end{proof}

\begin{lemma}  \label{annulus reduction}   Let $(M,T,\Delta)$ be an 
internal Mom-n structure on the compact hyperbolic manifold $N$.
Assume that every sheet of $\Delta$ is a disc and $\Delta$ is full.
If there exists an embedded annulus $A$ connecting the component $S$
of $\partial M-T$ to $T$, then there exists a non-elementary embedding
of a manifold $M'$ into $N$ such that $\partial M'$ is a union of tori
and $\rho_1(\Delta')+2\le \rho_1(\Delta)$.  Here $\Delta'$ is a handle
structure on $M'$ based on $T'$.  Either $T'=T$ or $T'$ is the union
of an essential annulus in $T$ and possibly a component of $\partial
M'$.
\end{lemma}

\begin{proof}  By Lemma \ref{reducible reduction} we  can
assume that $M$ is irreducible.  Since all sheets are discs and
letting $M_2=M$, we obtain from $\Delta$ a full handle structure $\Delta_2$
with no 0-handles based on the component $T_2=T$ of $\partial M_2$
where each 1 and 2-handle of $\Delta_2$ is of valence $\ge 3$.  Finally
$\rho_1(\Delta_2)\le \rho_1(\Delta)$.  We can assume that $A\subset
M_2$ is a least-weight normal annulus connecting $T_2$ to the boundary
component $S\neq T_2$.  Since $A$ is least weight, $A\cap T_2\times I$
is a union of discs and a single annulus.

Let $\partial A_0$ (resp.  $\partial A_1$) denote $A\cap T_2$
(resp. $A\cap S$).  Our $A$ has an induced handle structure $\Phi$
based on $\partial A_0$ as follows.  The base consists of the annular
component of $A\cap T_2\times [0,1]$, the 0-handles consist of the
disc components of $A\cap T_2\times[0,1]$, the 1-handles
(resp. 2-handles) consist of the intersections of $A$ with the
1-handles (resp. 2-handles).

Let $M_3$ denote the manifold obtained by splitting $M_2$ along $A$
and let $\Delta_3$ denote the induced handle structure.  As in \S 2,
the ball components of $T_2\times I$ split along $A$ to become 0-handles
of $\Delta_3$.  The remaining component is $B' \times I$, where $B'$
is $T_2$ split along $\partial A_0$ and hence $\Delta_3$ is based on
$B'\subset T_2$.  By Lemma \ref{pi image}, $M_3$ is non-elementarily
embedded in $N$.

If $\eta$ is a 1-handle of $\Delta_2$ and if $\{\eta_i\}$ denotes the
1-handles of $\Delta_3$ which descended from $\eta$, then $\sum_i
(\valence(\eta_i)-2) \le \valence(\eta)-2$ with equality if and only
if $A_1$ does not run over $\eta$.  In fact, counting with
multiplicity, if $A_1$ runs over the 1-handles of $\Delta_2$ more than
once then by operations as in Lemma \ref{free edge} we can pass to
$M_4$ and $\Delta_4$ with $\rho_1(\Delta_4)+2\le\rho_1(\Delta)$.  We
will now assume that $A_1$ ran over a unique 1-handle $\eta$ and it
did so with multiplicity one.  This implies that
$\rho_1(\Delta_3)+1=\rho_1(\Delta_2)$.  Also $A_1$ is the union of two
arcs $\alpha$ and $\beta$, where $\alpha$ lies in a 1-handle of
$\Delta_2$ and $\beta$ lies in a lake.  $\beta$ lies either in the
base of $\Phi$ or in a 0-handle $v^*$ of $\Phi$.

Give $A$ a transverse orientation. Call a 0-handle $v$ of $\Phi$
\emph{plus} (resp.  \emph{minus}) if the transverse orientation of the
disc $v\subset T_2\times I$ points away from (resp. towards) $T_2$.
Each such disc $v$ separates off, in $T_2\times I$, a 3-ball $v_B$.
Let $v_D$ denote $v_B\cap T_2\times 1$.

If $\beta$ lies in a 0-handle $v^*$, then Lemma \ref{0-handle
reduction} applies to $\Delta_3$ and the Lemma is proved.  Indeed,
since $A$ is least weight, the disc $v^*_D$ contains a bridge $b$ in
its interior and this bridge is not parallel to $\beta$.  In other
words, there
is no embedded disc in a lake whose boundary is a concatenation of four
arcs $\beta_1, \beta_2, \beta_3$ and $\beta_4$ such that $\beta_1=\beta,
\beta_3$ lies in $\partial b$, and $\beta_2$ and $\beta_4$ are arcs
lying in islands.  Otherwise, since $A_1$ runs over a unique 1-handle,
this implies that $\Delta_2$ had a valence-1 2-handle.  If $v^*_D$ is
disjoint from the other 0-handles of $\Phi$, then Lemma \ref{0-handle
reduction} applies to the 0-handle $v^*_B$ of $\Delta_3$.  Otherwise,
$v^*_B$ is split into balls by the various 0-handles of $\Phi$ and
Lemma \ref{0-handle reduction} applies to one of these balls.

 From now on we assume that $\beta$ lies in the base of $\Phi$.  Let
$X\subset T_2\times 1$ denote the union of the islands and bridges of
$\Delta_2$.  A similar but easier argument to the one given in the
previous paragraph shows that for each 0-handle $v$ of $\Phi$ either
$\partial v$ is boundary parallel in $X$ or Lemma \ref{0-handle
reduction} applies.  This implies that if $v\neq w$ are 0-handles of
$\Phi$, and $v_B\subset w_B$, then $\partial v$ and $\partial w$ are
normally parallel in $X$.  Furthermore, no 1-handle of $\Phi$ connects
a plus 0-handle to a minus 0-handle of $\Phi$.  Also, if $v^0$ and
$v^1$ are two 0-handles of $\Phi$ that are connected by a 1-handle,
then $v^0_B\cap v^1_B=\emptyset$.  Finally, there do not exist
0-handles $w^0, w^1, \cdots w^n$ of $\Phi$ such that for $i=1, 2,
\cdots, n-1,\ w^i$ is connected to $w^{i+1}$ by a 1-handle and
$w^n_B\subset w^0_B$.

It follows that there exists a disc $E\subset A$ whose boundary is the
union of two arcs $\phi$ and $\psi$ where $\phi$ is a proper arc in a
2-handle of $\Phi$ and $\psi$ lies in the $T\times 1\cap (\base(A))$.
Furthermore $E\cap (0,1\text{-handles}\ \Phi)\neq \emptyset$ and is
connected.  By the previous paragraph if $v,w$ are 0-handles of $\Phi$
lying in $E$, then they are of the same parity and $v_B\cap
w_B=\emptyset$.

Therefore, $E$ can be normally isotoped, with respect to $\Delta_2$ to
a disc $G\subset \partial M$ such that $\partial G$ is a union of two
arcs, one lying in a lake and the other in a 2-handle.  Now apply (ii)
of Lemma \ref{boundary reduction} to $\Delta_2$.\end{proof}

\begin{lemma} \label{annulus reduction II}  Let $(M,T,\Delta)$ be a 
full internal Mom-$n$ structure on the compact hyperbolic 3-manifold $N$
such that every sheet is a disc.  Suppose that there exists an
essential embedded annulus $A$ with $\partial A\cap T=\emptyset$.
Then either there exists a full internal Mom-$k$ structure
$(M',T',\Delta')$ on $N$ with $k<n$ or there exists a non-elementary embedding
$M'\to N$ with handle structure $(M',T',\Delta')$ where $\partial M'$
is a union of tori, $T'$ is a component of $\partial M'$, and
$\rho_1(\Delta')+2\le \rho_1(\Delta)$.\end{lemma}

\begin{proof}  Assume that $M$ is irreducible.  By Lemma \ref{annulus 
reduction} we can assume that no essential annulus connects $T$ to a
component of $\partial M-T$.  Let $A$ be a least-weight essential
normal annulus with $A\cap T=\emptyset$.  By the second sentence,
$A\cap T\times I$ is a union of discs.  Let $M_1$ be the non-elementary
component of M split along $A$ and let $(M_1,T_1,\Delta_1)$ be the induced
handle structure.  As in the previous proof, $\rho_1(\Delta_1)+1\le
\rho_1(\Delta)$.  Indeed, since $\partial A$ traverses at least two
1-handles of $\Delta$, counted with multiplicity, the inequality will
be strict unless $A$ normally double covers a Mobius band.  If $M_1$
is disjoint from $T$, then the result follows from Lemma \ref{0-handle
reduction}.  Assume now that $T\subset\partial M_1$ and $\Delta_1$ is
full.

Cancel the 0-handles with 1-handles to obtain $(M_2,T_2,\Delta_2)$.
If $\Delta_2$ is not full, then some 1-handle of valence $\ge 3$ was
cancelled and hence $\rho_1$ is reduced.  Next cancel valence-1
2-handles to obtain $(M_3,T_3,\Delta_3)$.  Again fullness is preserved
or $\rho_1$ is reduced.  If $\Delta_3$ has a valence-2 2-handle
$\sigma$ connecting distinct 1-handles, then cancel one of these
1-handles with $\sigma$.  If the 2-handle $\sigma$ attaches to the
same 1-handle $\alpha$, then the core of $\sigma\cap\alpha$ can be
viewed as an embedded annulus or Mobius band $C$ in $M_3$.  Splitting
$M_3$ along $C$ we obtain $(M_4,T_4,\Delta_4)$ where $T_4$ is the
newly created boundary component of $M_4$ and $\Delta_4$ is the
induced structure.  Note that $\rho_1(\Delta_4)<\rho_1(\Delta_3)$
since $\Delta_3$ is full.  After finitely many such operations the
lemma is proved.  \end{proof}

\begin{lemma}\label{annulus reduction III}  Let $(M,T,\Delta)$ be a 
full internal Mom-$n$ structure on the compact hyperbolic 3-manifold $N$
such that every sheet is a disc.  Suppose that there exists an
essential embedded annulus $A$ with $\partial A\subset T$.  Then
either there exists a full internal Mom-$k$ structure
$(M',T',\Delta')$ on $N$ with $k<n$ or there exists a non-elementary embedding
$M'\to N$ with handle structure $(M',T',\Delta')$ where $\partial M'$
is a union of tori, $T'$ is a component of $\partial M'$, and
$\rho_1(\Delta')+2\le \rho_1(\Delta)$.\end{lemma}

\begin{proof}  If $M$ is reducible or there exists an essential 
annulus with some boundary component disjoint from $T$, then apply
Lemma \ref{reducible reduction}, Lemma \ref{annulus reduction} or
Lemma \ref{annulus reduction II}.  Let $(M,\partial M-T, \Sigma)$ be
the dual handle structure.  I.e. the 1-handles (resp. 2-handles) of
$\Sigma$ are in 1-1 correspondence with the 2-handles
(resp. 1-handles) of $\Delta$ and for each handle $\sigma$, the core of
$\sigma$ is the cocore of the dual handle and vice versa.  Note that
$\rho_1(\Sigma)=\rho_1(\Delta)$ and each 1-handle of $\Sigma$ is of
valence $3$, but $\Sigma$ may have valence-2 2-handles.

Let $A$ be an essential annulus, least weight with respect to
$\Sigma$.  Since each essential annulus has its entire boundary in
$T$, $A\cap (\partial M-T)\times I$ is a union of discs.  Let $M_1$ be
the non-elementary component of $M$ split along $A$ with
$(M_1,R_1,\Sigma_1)$ the induced structure.  As before
$\rho_1(\Sigma_1)<\rho_1(\Sigma)$.  Eliminate the 0-handles and low
valence-1 and 2-handles of $\Sigma_1$ to obtain $(M_2,R_2,\Sigma_2)$
and note that either $\Sigma_2$ is full or
$\rho_1(\Sigma_2)+2\le\rho_1(\Sigma)=\rho_1(\Delta)$.  If $\Sigma_2$
is full, then $(M_3,T_3,\Delta_3)$ the handle structure dual to
$(M_2,R_2,\Sigma_2)$ is full.  Let $(M_4,T_4,\Delta_4)$ be obtained by
eliminating the valence-2 2-handles as in the previous proof.  As
before $\rho_1(\Delta_4)+1\le \rho_1(\Delta)$ with equality holding
only if $\Delta_4$ is full.
\end{proof}

\begin{lemma} \label{torus reduction}  Let $(M,T,\Delta)$ be a full 
internal Mom-$n$ structure on the compact hyperbolic 3-manifold $N$.  If
$M$ is not hyperbolic, then either there exists a full internal Mom-$k$
structure $(M',T',\Delta')$ on $N$ with $k<n$ and $M'$ hyperbolic or
there exists a general based internal Mom-l, $l+2\le n$, structure
$(M',B',\Delta')$ on $N$.\end{lemma}

\begin{proof}  By Lemma
\ref{reducible reduction}, \ref{sheet reduction}, \ref{annulus
reduction}, \ref{annulus reduction II} or \ref{annulus reduction III}
it suffices to consider the case that $M$ is irreducible and
anannular.  By Thurston \cite{Th2}, if $M$ is not hyperbolic, then it
contains an essential torus.  Let $F$ be a least-weight essential
torus.  Since there are no essential annuli, $F\cap T\times I$ is a
union of discs.

  Let $M_1$ denote the component of $M$ split along $F$ which contains
$T$ and let $M_1'$ denote the other component. Let $\Delta_1$ and
$\Delta_1'$ denote the induced handle structures.  Note that
$\Delta_1$ is based on $T$, $\Delta_1'$ is a classical structure and
$\rho_1(\Delta_1)+\rho_1(\Delta_1')=\rho_1(\Delta)$.  By Lemma \ref{pi
image} one of $M_1$ or $M_1'$ is non-elementary .  Let $X\subset
T\times 1$ be the union of the islands and bridges of $\Delta$,
$Y=M_1'\cap X$ and $Z=M_1\cap X$.

If $M_1$ is elementary, $\rho_1(\Delta_1)\ge 1$ since $Z$ is
nonplanar.  In the usual way obtain a handle structure $\Sigma$ on
$M_1'$ such that $\rho_1(\Sigma)\le \rho_1(\Delta')$, $\Sigma$ has
exactly one 0-handle and each 1-handle is of valence $\ge 2$.  Since
$M_1'$ is non-elementary, some 1-handle is of valence $\ge 3$.  To
complete the proof apply Lemma \ref{0-handle reduction}.

 From now on assume that $M_1$ is non-elementary.  We first show that
$\rho_1(\Delta_1')\ge 1$.  If $\rho_1(\Delta_1')=0$, then $Y$ is a
union of annuli, each component of which intersects $\partial X$ in
$\le 1$ circle, since $\Delta$ is full.  If $Y\cap \partial
X\neq\emptyset$, then each component must have this property and hence
$F$ is boundary parallel in $M$.  If $Y\cap \partial X=\emptyset$,
then $M_1'$ is an I-bundle over a Klein bottle, which is impossible in
an orientable hyperbolic 3-manifold.

In the usual way pass from $(M_1,T_1,\Delta_1)$ to a Mom structure
$(M_2,T_2,\Delta_2)$ with $\rho_1(\Delta_2)\le
\rho_1(\Delta_1)$. Since $M$ is anannular, $\Delta_2$ has no annular
lakes.  If some component of $\partial M_2$ is a convolutube, then
reimbed $M_2$ in $N$ to get a full internal Mom-k structure where
$k<n$.\end{proof}

\begin{section} {From  Mom-$n$ to  Hyperbolic
Mom-$n$}\end{section}

The following is the main result of this section:

\smallskip
\begin{theorem} \label{hyperbolic Mom-n} If $(M,T,\Delta)$ is a full 
internal Mom-$n$ structure on the compact hyperbolic 3-manifold $N$ and
$n\le 4$, then there exists a full internal Mom-$k$ structure
$(M',T',\Delta')$ on $N$ where $M'$ is hyperbolic and $k\le n$.

For general $n$, either there exists an internal Mom-$k$ structure
$(M',T', \Delta')$ in $N$ with $k\le n$ and $M$ hyperbolic or there
exists a general based internal Mom-$k$ structure $(M',B',\Delta')$
on $N$ such
that $M'$ is hyperbolic and each component of $B'$ is either a component
of $\partial M'$ or an annulus which is essential in a component of $\partial
M$.  If $M\neq M'$, then $k+2\le n$.\end{theorem}

We first prove some preliminary lemmas about complexity (see Definition \ref{def:valence}).

\begin{lemma}(Clean-Up Lemma)\label{clean up}  Let $M$ be a compact 
3-manifold with $\partial M$ a union of tori and $N$ a compact
hyperbolic 3-manifold.  If $M\to N$ is a non-elementary embedding and
$\Delta$ is a handle structure on $(M,B)$, then there exists
$(M',B',\Delta')$ a weak general based internal Mom-$n$ structure on $N$
with $C(\Delta')\le C(\Delta)$.

If $M$ is hyperbolic, then either $\rho_1(\Delta')<\rho_1(\Delta)$ or
$M'=M$.
\end{lemma}

\begin{proof} Proof by induction on $C(\Delta)$.  By the usual
cancellation operations pass to a handle structure $\Delta_2$ on
$(M,B)$ without 0-handles or 1 and 2-handles of valence 1.  If $\Delta$
has a valence-0 1-handle, then delete it and proceed as in Lemma
\ref{free edge} to obtain $\Delta_3$ on $(M_3,B_3)$ with $M_3$
non-elementary and $\partial M_3$ a union of tori.  All of these
operations reduce $C(\Delta)$.

If $\Delta_2$ has a valence-0 2-handle $\sigma$, then compress $M_3$
along a disc $D$ which passes once through $\sigma$ and intersects
$B\times I$ in an annulus.  Let $\Delta_4$ be the resulting handle
structure on $(M_4,B_4)$ where $M_4$ is the non-elementary component of
$M_3$ split along $D$ and $B_4$ is $B_3\cap M_4$.  This splits
$\sigma$ into two 2-handle components, at least one of which lies in
$\Delta_4$.  Absorb such 2-handles into $B_4\times I$.  This creates
$(M_5,B_5,\Delta_5)$ where $M_5=M_4$ and $B_5$ is obtained by
attaching one or two discs to $\partial B_4$.  If some
component of $B_5$ is a 2-sphere $S$, then create $(M_6,B_6,\Delta_6)$
by filling in $S$ with a 3-cell $E$ and identifying $S\times I\cup E$
as a 0-handle of $\Delta_6$.  Now cancel with a 1-handle to obtain
$(M_7,B_7,\Delta_7)$ with $C(\Delta_7)< C(\Delta_3)$.  If one or two
components of $B_5$ are discs and there exists a non-simply connected
component of $B_5$, then the ball components of $B_5\times I$ are now
re-identified as 0-handles and cancelled with 1-handles to obtain
$(M_7,B_7,\Delta_7)$.  If all the components of $B_5$ are discs, then
transforming the components of $B_5\times I$ to 0-handles produces a
classical handle structure $\Delta_{5.5}$ on $M_5$.  Now apply the
proof of Lemma \ref{0-handle reduction} to create $(M_6,B_6,\Delta_6)$
with $C(\Delta_6)$ reduced.  If some component of $\partial M_6$ is a
2-sphere, then delete a 2-handle facing it to obtain
$(M_7,B_7,\Delta_7)$.

By repeatedly applying the above operations we can assume that our
$\Delta_7$ has no 0-handles, and each 1-handle and 2-handle is of
valence $\ge 2$.  Furthermore, no component of $B_7$ is a disc.  If
some component $G$ of $B_7$ is not $\pi_1$-injective in $\partial
M_7$ then attach a 2-handle to a component of $\partial G\times I$
whose restriction to $\partial M_7$ bounds a disc, to create
$(M_8,B_8,\Delta_8)$.  Delete a 2-handle which faces the resulting
2-sphere boundary component and simplify as in the previous paragraphs
to create $(M_9,B_9,\Delta_9)$ with $C(\Delta_9)<C(\Delta_7)$.  We can
now assume that each component of $B_9$ is essential in $\partial
M_9$, $\Delta_9$ is 0-handle free and all 1 and 2-handles are valence
$\ge 2$.

  If $\partial M_9$ contains convolutubes, then reimbed $M_9$ in $N$
to eliminate them.  If $M$ is not hyperbolic, then the resulting
$(M_{10},B_{10},\Delta_{10})$ satisfies the conclusions of our clean
up lemma.

If $M$ is hyperbolic, then the proof also follows as above.  However,
note that $\partial M$ is incompressible so in the above process
valence-0 1-handles or lakes compressible in $B\times I$ can be
eliminated without changing $M$ and with reducing complexity.  The
topology of $M$ may change if we delete a 2-handle or apply Lemma
\ref{0-handle reduction}; however, if such operations must be done,
they can be done in such a manner that reduces $\rho_1$.
\end{proof}

\begin{lemma}  \label{weak} If
$(M,B,\Delta)$ is a strictly weak general based Mom-$n$ structure with
$M$ hyperbolic, then there exists a general based Mom-$k$
structure $(M,B',\Delta')$ with $C(\Delta')<C(\Delta)$.
\end{lemma}

\begin{proof}  Proof
by induction on $C(\Delta)$.  If $\sigma$ is a valence-2 2-handle
which goes over distinct 1-handles, then cancelling a 1-handle with
$\sigma$ creates $(M,B, \Delta_1)$ with $C(\Delta_1)<C(\Delta)$.

We now assume that no valence-2 2-handle $\sigma$ goes over distinct
1-handles.  Suppose that $\sigma$ goes over the same 1-handle $\alpha$
twice.  Then $\sigma\cup\alpha$ can be viewed as an embedded annulus
or Mobius band $A$ with boundary on $B$.  Since $M$ is hyperbolic and
orientable, $A$ is either a boundary parallel annulus or a
compressible annulus.

If $A$ is boundary parallel, then it together with an annulus on $T$
cobound a solid torus $V\subset B$ such that $A$ wraps longitudinally
around $V$ exactly once.  Melting $V$ into $B\times [0,1]$ eliminates
both $\alpha$ and $\sigma$, together with all the 2-handles and
1-handles inside of $V$.  The resulting manifold $(M,B_1,\Delta_1)$ is
a weak provisional Mom-$k$ with $k\le n$ whose induced handle structure
$\Delta_1$ satisfies $C(\Delta_1)<C(\Delta)$.

If $A$ is not boundary parallel, then $M$ is reducible, a
contradiction.
\end{proof}

\begin{lemma}  \label{lower bound} If $(M,B,\Delta)$ is a general 
based Mom-$n$ and $M$ is hyperbolic, then $n\ge 2$.
\end{lemma}

\begin{proof}  This  follows by direct calculation.\end{proof}

\begin{lemma}\label{general base 2}  If the compact hyperbolic 
3-manifold $N$ has a general based internal Mom-2 structure
$(M,B,\Delta)$, then it has a full internal Mom-2
structure.\end{lemma}

\begin{proof}  By splitting along annuli in $B\times I$ we can assume 
that every non-peripheral lake of $B$ is a disc.  If $B$ is the union
of two tori $T_1$ and $T_2$, then each 1-handle must connect $T_1$ to
$T_2$.  This implies that each 2-handle is of even valence which is a
contradiction.  If $B$ consists of a single torus, then $(M,B,\Delta)$
is a full Mom-2 structure.  We finally assume that $B$ contains an
annulus.  We only discuss the case that $B$ is connected; the case
where $B$ is either the union of an annulus and torus or two annuli is
similar and easier.

First suppose that $\Delta$ has a single sheet of valence four.  Let
$\eta$ denote the valence-4 1-handle of $\Delta$. There exists an
essential compressing disc $D$ for $B\times I$ which cuts across the
bridges in at most three components.  View $N(D)$ as a 1-handle and
$B\times I-\cirN(D)$ as a 0-handle to obtain a classical handle
structure with two 1-handles respectively of valence 4 and $\le 3$.
Now as in Lemma \ref{0-handle reduction} hollow out the 0-handle and
$\eta$ to get a non-elementary $M_1$ with handle structure $\Delta_1$
based on a torus with a single 1-handle of valence $\le 3$.  Therefore
we obtain a Mom-1 structure on a hyperbolic 3-manifold, which is a
contradiction.

Now consider the case that $\Delta$ has two valence-3 1-handles.  If
some essential disc $D$ in $B\times I$ cuts the bridges in $\le 2$
components, then as above, we obtain a handle structure with one
0-handle and three 1-handles.  Hollowing out a valence-3 1-handle and
the 0-handle produces a handle structure $\Delta_1$ on a
non-elementary manifold $M_1$ with $\rho_1(\Delta_1)=1$, which is a
contradiction.  If no essential disc $D$ exists as above, then we can
find one which cuts the bridges in exactly three components.  One
readily enumerates the possible handle structures that satisfy our
assumptions of valence-3 1-handles and non-peripheral disc lakes.
After applying the hollowing out procedure, a Mom-2 handle structure
is created.  Of the two choices of which 1-handle is to be hollowed
out, one will produce a full Mom-2 structure.

\end{proof}

\begin{remark} As an example of the type of non-full Mom-2 discussed
  in the previous lemma, note that the figure-8 knot complement
  contains a non-full Mom-2 structure with two 1-handles of valence 3,
  in addition to having a full Mom-2 structure.
\end{remark}

\noindent\emph{Proof of Theorem \ref{hyperbolic Mom-n}} The proof is
by induction on $\rho_1(\Delta)$.  If $(M,T,\Delta)$ is not hyperbolic,
then apply Lemma \ref{torus reduction} to obtain $(M_1,T_1,\Delta_1)$
where $M_1\to N$ is a non-elementary embedding, $\partial M_1$ is a
union of tori and either $(M_1,T_1,\Delta_1)$ is a full internal Mom-$k$
structure, $k<n$, or $\rho_1(\Delta_1)+2\le \rho_1(\Delta)$. In the former
case the proof follows by induction. In the latter case as seen in the next
paragraph we will produce a general based internal Mom-$k$ structure
$(M_4,B_4,\Delta_4)$ with $k\le n-2$ and $M_4$ hyperbolic.  Lemma
\ref{lower bound} then implies that $k\ge 2$ and hence $n\ge 4$.  If
$k=2$, then Lemma \ref{general base 2} implies that $N$ has a full
internal Mom-2 structure $(M_5,B_5,\Delta_5)$.  $M_5$ must be
hyperbolic, or else the above arguement for $n=2$ will give a
contradiction to Lemma \ref{lower bound}, completing the proof.

So assume that $(M_1,T_1,\Delta_1)$ is not a full internal Mom
structure.  If $M_1$ is reducible, then split along a normal reducing
2-sphere, retain the non-elementarily embedded component, cap off the
resulting 2-sphere boundary component with a 3-handle and cancel that
3-handle with a 2-handle.  After a sequence of such operations we
obtain $(M_2,T_2,\Delta_2)$ with $\rho_1(\Delta_2)\le\rho_1(\Delta_1)$
and $M_2$ irreducible.  If $M_2$ contains an embedded essential torus
$R$, then split along $R$ and retain the non-elementarily embedded
component.  After a sequence of such operations we obtain
$(M_3,T_3,\Delta_3)$ with $\rho_1(\Delta_3)\le \rho_1(\Delta_2)$ and
$M_3$ is irreducible and geometrically atoroidal.  Since $M_3$ is
non-elementarily embedded in $N$, it is not a Seifert fibered space and
hence is anannular and so by Thurston \cite{Th2} it is hyperbolic.
After repeatedly applying Lemmas \ref{clean up} and \ref{weak} we
obtain a general based internal Mom-k structure $(M_4,T_4,\Delta_4)$
on $N$ where $M_4$ is hyperbolic and $\rho_1(\Delta_4)\le
\rho_1(\Delta_3)$. \qed

\end{section}

\begin{section}{Enumeration of hyperbolic Mom-$n$'s
for $2\le n \le 4$}\label{enumeration}
\end{section}

Let \((M,T,\Delta)\) be a full hyperbolic Mom-$n$, with \(2\le n\le
4\). The handle
structure $\Delta$ collapses to a cellular complex $K$ in the
following fashion. Each 1-handle
collapses to the arc at its core, and each 2-handle collapses to the
disc at its core (expanded as necessary so that it is still attached
to the cores of the appropriate 1-handles). Also, $T\times I$
collapses to $T\times 1$, subdivided into 0-cells, 1-cells, and
2-cells corresponding to the islands, bridges, and lakes of
$(M,T,\Delta)$. (Note that if $(M,T,\Delta)$ were not full, we might
have a non-simply connected lake and $K$ would not be a proper
cellular complex.)

The resulting complex $K$ is a \emph{spine} for $M$ in the sense of
\cite{MF}. If all of the 1-handles of $\Delta$ are of valence 3, then
it is also a \emph{special spine} in the sense of \cite{MF}; however
$K$ is not a special spine in general. In particular, in a special
spine the link of each point is either a circle or a circle with two
or three radii, but if $\Delta$ has a 1-handle of valence $n$ then the
endpoints of the corresponding arc in $K$ will have links which are a
circle with $n$ radii. This, however, is the only way in which $K$
fails to be a special spine.

In section 2 of \cite{MF}
Matveev and Fomenko describe how a manifold with a special spine can
be reconstructed by gluing together truncated or ideal simplices dual
to the vertices of the spine. This construction is easily generalized
to our situation, and shows that $M$ can be reconstructed from $K$ by
gluing together ideal polyhedra dual to the vertices of $K$. The
result is an ideal cellulation of $M$ which is dual to the cellular
complex $K$.


The $3$-cells of this cellulation will be dual to the elements of
\(K^0\), which consist of the endpoints of the cores of the
$1$-handles of $\Delta$.  In addition, since we've assumed each
$1$-handle of $\Delta$ meets at least two $2$-handles, each point
\(v\in K^0\) will be the endpoint of at least two curves in \(T\times
1 \cap K^1\). Hence if \(n_v\) is the valence of $v$ in the 1-skeleton
of $K$ then
\(n_v\ge 3\).  If \(n_v \ge 4\) then $v$ is dual to an
\((n_v-1)\)-sided pyramid: the base of the pyramid is dual to the core
of a $1$-handle while the sides are dual to curves in \(T\times 1 \cap
K^1\). If \(n_v = 3\) then $v$ is dual to a ``digonal pyramid'', which
we eliminate from the cellulation by collapsing it to a face in the
obvious fashion.
Thus $K$ is dual to a cellulation of $M$ by ideal pyramids. Since the
bases of these pyramids correspond to the ends of the $1$-handles of
$\Delta$, we can pair them up into a collection of ideal dipyramids.


We can say more concerning the possible types and combinations of
dipyramids. On one hand, each vertex $v$ is adjacent to \(n_v-1\)
edges in \(T\times 1\cap K^1\), and each such edge has two endpoints;
on the other hand, the core of each $2$-cell of $\Delta$ contributes
three edges to \(T\times 1 \cap K^1\), and there are $n$ such cores in
a Mom-$n$.  Therefore \(\Sigma_v (n_v-1)=6n\) in a
Mom-$n$. Furthermore, \(n_v-1\) must be at least $2$ and (if it's
greater than $2$) equals the number of sides of the pyramid dual to
$v$. Finally the vertices $v$ occur in pairs since each one
corresponds to an end of a $1$-handle, and the vertices in each pair
have the same valence. Therefore for a Mom-2, there are only two
possibilities: four three-sided pyramids, which glue together to form
two three-sided dipyramids, or two four-sided pyramids and two
``digonal pyramids'', which (after eliminating the ``digonal
pyramids'') glue together to form a single ideal
octahedron. Similarly, there are only three possibilities for a Mom-3:
three three-sided dipyramids, a three-sided dipyramid together with an
octahedron, or a five-sided dipyramid by itself. The five
possibilities for a Mom-4 are: four three-sided dipyramids, two
three-sided dipyramids and an octahedron, one three-sided dipyramid
and one five-sided dipyramid, two octahedra, or one six-sided
dipyramid.

Thus, if $(M,T,\Delta)$ is a hyperbolic Mom-2, Mom-3, or Mom-4 then
$M$ can be obtained by gluing together the faces of one of these ten
sets of ideal polyhedra. Enumerating the possibilities for $M$ then
becomes a matter of enumerating the ways in which the faces of these
polyhedra can be glued together to form a hyperbolic $3$-manifold.

This task is simplified somewhat by the following observation: the
faces of each dipyramid always have exactly one vertex which is dual
to the cusp neighborhood \(T\times [0,1)\).  When gluing the polyhedra
together to form $N$, all such vertices must be identified with one
another and with no other vertices. Thus given any two faces, there is
only one possible orientation-preserving way that those two faces
could be glued together.


Hence it is sufficient to enumerate the number of ways in which the
faces of one of the ten sets of polyhedra can be identified in
pairs. Although it is almost trivial to program a computer to do this,
care must be taken as the number of possibilities is a factorial
function of the number of faces, and a naive approach can rapidly
exhaust a computer's memory. To reduce the demands on the computer, a
refinement to the naive approach was employed. First, for each
possible set of polyhedra a symmetry group was computed. Each
dipyramid has dihedral symmetry, while if a given set of polyhedra
contains two dipyramids with the same number of sides then they can be
exchanged to provide an additional symmetry. Secondly, an ordering was
chosen for the set of all possible pairings of faces, namely the
lexicographic ordering of the pairings when represented as
permutations. Our computer program considered the set of pairings in
order, and any pairing was immediately rejected if it was conjugate to
a previous pairing via an element of the symmetry group. This
considerably reduced the running time of the program.


The next step in the process is to eliminate pairings which result in
obviously
non-hyperbolic manifolds.  While the program SnapPea can in principle
handle this, for reasons of speed our program checked one necessary
criterion itself: whether the link of every ideal vertex was
Euclidean. Computing the Euler characteristic of the link of each
ideal vertex in the cellular complex resulting from a pairing was easy
to do and eliminated many cases from consideration. Our program also
eliminated any pairing in which the vertices supposedly dual to the
original cusp neighborhood or solid torus in fact glued together to
form two or more ideal vertices.

The above considerations resulted in a list of gluing descriptions of
4,231 manifolds which might be hyperbolic Mom-2's or Mom-3's.
At this point, SnapPea was employed to try and compute
hyperbolic metrics for each of these manifolds, and to find further
hyperbolic symmetries among the manifolds which admitted such
metrics. SnapPea claimed to find hyperbolic metrics in 164 cases. In
three of those cases SnapPea had experienced an obvious floating-point
error and ``found'' a hyperbolic metric with an absurdly low
volume. In those three cases the programs Regina and GAP (which,
unlike SnapPea, do not rely on floating-point arithmetic) were used to
calculate the fundamental groups of the corresponding manifolds. In
all three cases the fundamental group was isomorphic to the group
\(\langle a,b|[a,b^n]\rangle\) where \(n=3\) or $5$, which has a
non-trivial center. Therefore, these three manifolds could not
possibly be hyperbolic and were rejected. That left 161
cases, which SnapPea identified as belonging to a total of 21 isometry
classes of hyperbolic manifolds.

Some comments about rigor are in order here. Since SnapPea relies on
floating-point arithmetic, some of its results are unavoidably
inexact. In particular, there is no guarantee that SnapPea will find
a hyperbolic metric on a manifold even if one exists, or that SnapPea
will correctly discern the absence of a hyperbolic metric in cases
where it doesn't
exist. In practice it is our experience that if one is careful to
allow SnapPea to simplify a triangulation before attempting to find a
metric, then if a metric exists SnapPea will either find it or fail to
make a determination, while if a metric doesn't exist SnapPea will
either correctly say so or on rare occasions ``find'' a metric with
absurdly low volume due to floating-point error. Still, from a
standpoint of rigor this is problematic. Fortunately there is at
least one task which SnapPea does perform exactly, and that is
finding isometries between two different cusped manifolds: SnapPea will only
report that an isometry exists if it finds identical triangulations of
the two manifolds. (See \cite{We}; in particular see the comments in
the source code file \emph{isometry.c}.) This is a combinatorial operation,
not a floating-point one, and hence we are confident that SnapPea
performs this operation rigorously.

Those familiar with SnapPea's source code may object that SnapPea
re-triangu\-lates each manifold before determining if an isometry
exists, and that SnapPea uses floating-point information to choose the
re-triangu\-lation. To this objection we would reply that while
floating-point information is used to \emph{choose} the
re-triangulation, the actual re-triangu\-lating is still a combinatorial
operation, i.e. it uses integer arithmetic. The new triangulation is
guaranteed to have the same topological type (see the comments in
\emph{canonize\_part\_1.c} from \cite{We}), and hence the
possibility of floating-point error does not invalidate the result
when SnapPea reports that it has found an isometry.

Thus while we are trusting SnapPea when it says that the 161 manifolds
mentioned above are all in the isometry class of one of 21 manifolds
from the SnapPea census, we are confident that we are not sacrificing
rigor in so doing. Furthermore, the census manifolds were
recently confirmed to be hyperbolic by Harriet Moser in \cite{Mos},
establishing that we have found 21 different hyperbolic Mom-2's and
Mom-3's.

Unfortunately, we still can't trust SnapPea when it fails to
find an hyperbolic metric for a given manifold,
as that result is not guaranteed to be
rigorous. This means that there are 4,067 manifolds from the above
list of 4,231 which may still be hyperbolic despite SnapPea evidence
to the contrary. These manifolds were analyzed separately in the same
way as the three manifolds for which SnapPea claimed to have found an
absurdly low hyperbolic volume. Namely,
we used Regina and GAP as before to compute the fundamental
groups of the manifolds in question, and then examined the list of groups
to see if any of them might be the fundamental group of a hyperbolic
manifold. The vast majority of the groups on the list either had a
non-trivial center, or else had two rank-2 Abelian subgroups
which intersected in a rank-1 Abelian subgroup (also impossible in the
fundamental group of a hyperbolic 3-manifold). Some of the groups required
further analysis but were still eventually rejected; for example, many
groups had an index-two subgroup with one of the above properties even
when it was not clear that the whole group had such properties.

In the end the hand analysis did not reveal any new hyperbolic 3-manifolds
in the list of gluing descriptions. This completes the proof of the
following:

\smallskip
\begin{theorem} \label{list of mom-2's and mom-3's}
There are 3 hyperbolic manifolds $M$ such that \((M,T,\Delta)\) is a
Mom-2 for some $T$ and $\Delta$: the manifolds known in SnapPea's
notation as m125, m129, and m203. There are 18 additional hyperbolic manifolds
$M$ such that \((M,T,\Delta)\) is a Mom-3 structure for some $T$ and
$\Delta$: the manifolds known in SnapPea's notation as m202, m292,
m295, m328, m329, m359, m366, m367, m391, m412, s596, s647, s774,
s776, s780, s785, s898, and s959.
\end{theorem}
Some comments about this list are in order. The manifold m129, better
known as the complement of the Whitehead link, is the only manifold on
this list which is obtained by gluing together the faces of an ideal
octahedron. Also, all but one of these manifolds have two cusps. The
exception is the three-cusped s776, which is the complement in \(S^3\)
of a three-element chain of circles (the link \(6_1^3\) in Rolfsen's
notation).

Enumerating hyperbolic Mom-4's was more difficult: merely
enumerating the possible gluing descriptions resulted in a list of 1,033,610
possibilities (compared to 4,231 possibilities in the previous case).
From this list, SnapPea identified 138 different
hyperbolic manifolds. In another 493 cases, SnapPea was either unable
to make a determination or else experienced an obvious floating-point
error. In each such case, the fundamental group of the corresponding
manifold was again computed by Regina and GAP, and in each case the
fundamental group was isomorphic to \(\langle a,b |[a^n,b]\rangle\)
for some \(n\), or else had two rank-2 Abelian subgroups which
intersected in a rank-1 Abelian subgroup.  Therefore these exceptional
cases do not correspond to hyperbolic manifolds. Note that all of the
Mom-2's and Mom-3's appear in the Mom-4 list; the same manifold can
admit multiple handle structures.

Based on the above result, we propose the following:

\begin{conjecture}\label{number of mom-4's}
There are 138 hyperbolic manifolds $M$ such that \((M,T,\Delta)\) is a
Mom-2, Mom-3, or Mom-4 for some $T$ and $\Delta$. Of these, 117 are
\emph{strict} Mom-4's, i.e. Mom-4's which are not Mom-2's or Mom-3's.
\end{conjecture}

Of the 117 strict Mom-4's, SnapPea was successfully used to identify
83 of them as manifolds from the SnapPea census. Those manifolds appear
in
\begin{figure}
\begin{center}
\begin{tabular}{|c|c|c|c|c|c|c|}
\hline\hline
m357 & s579 & s883 & v2124 & v2943 & v3292 & v3450 \\
m388 & s602 & s887 & v2208 & v2945 & v3294 & v3456 \\
s441 & s621 & s895 & v2531 & v3039 & v3376 & v3468 \\
s443 & s622 & s906 & v2533 & v3108 & v3379 & v3497 \\
s503 & s638 & s910 & v2644 & v3127 & v3380 & v3501 \\
s506 & s661 & s913 & v2648 & v3140 & v3383 & v3506 \\
s549 & s782 & s914 & v2652 & v3211 & v3384 & v3507 \\
s568 & s831 & s930 & v2731 & v3222 & v3385 & v3518 \\
s569 & s843 & s937 & v2732 & v3223 & v3393 & v3527 \\
s576 & s859 & s940 & v2788 & v3224 & v3396 & v3544 \\
s577 & s864 & s941 & v2892 & v3225 & v3426 & v3546 \\
s578 & s880 & s948 & v2942 & v3227 & v3429 &  \\
\hline\hline
\end{tabular}
\end{center}
\caption{Conjectured list of SnapPea manifolds which are strict Mom-4's.}
\end{figure}
Figure 1. SnapPea was not able to identifiy the remaining 34
manifolds, and in fact 33 of those manifolds have volumes which do not
appear anywhere in the SnapPea census, presumably because the Matveev
complexity of the corresponding manifolds is greater than 7 (see
\cite{MF}). The remaining manifold has the same volume and homology as
the census manifold v3527; it is conceivable that SnapPea was simply
unable to find a corresponding isometry.

The unidentified manifolds are listed in
\begin{figure}
\begin{center}
\begin{tabular}{|c|}
\hline\hline
$(3,3,4\ ;\ 3,6,8,0,13,19,1,15,2,17,14,18,16,4,10,7,12,9,11,5)$ \\
$(3,3,4\ ;\ 3,6,11,0,10,9,1,15,14,5,4,2,16,18,8,7,12,19,13,17)$ \\
$(3,3,4\ ;\ 3,6,11,0,9,18,1,19,13,4,15,2,16,8,17,10,12,14,5,7)$ \\
$(3,3,3,3\ ;\ 3,6,9,0,13,19,1,22,14,2,17,12,11,4,8,23,18,10,16,5,21,20,7,15)$ \\
$(3,3,3,3\ ;\ 3,4,6,0,1,19,2,9,13,7,14,15,23,8,10,11,20,21,22,5,16,17,18,12)$ \\
$(3,3,3,3\ ;\ 3,6,12,0,17,9,1,11,18,5,23,7,2,20,15,14,21,4,8,22,13,16,19,10)$ \\
$(3,3,4\ ;\ 3,4,6,0,1,18,2,11,16,10,9,7,15,19,17,12,8,14,5,13)$ \\
$(3,3,4\ ;\ 3,4,6,0,1,18,2,14,12,13,19,17,8,9,7,16,15,11,5,10)$ \\
$(3,3,4\ ;\ 3,4,6,0,1,18,2,14,16,13,19,17,15,9,7,12,8,11,5,10)$ \\
$(3,3,4\ ;\ 3,6,10,0,13,9,1,16,19,5,2,17,14,4,12,18,7,11,15,8)$ \\
$(3,3,4\ ;\ 3,6,10,0,8,13,1,16,4,18,2,17,15,5,19,12,7,11,9,14)$ \\
$(3,3,4\ ;\ 3,6,10,0,8,13,1,18,4,16,2,19,15,5,17,12,9,14,7,11)$ \\
$(4,4\ ;\ 15,10,13,8,11,14,9,12,3,6,1,4,7,2,5,0)$ \\
$(4,4\ ;\ 15,14,5,6,9,2,3,10,13,4,7,12,11,8,1,0)$ \\
$(4,4\ ;\ 15,14,9,8,11,10,13,12,3,2,5,4,7,6,1,0)$ \\
$(4,4\ ;\ 15,4,13,6,1,8,3,10,5,14,7,12,11,2,9,0)$ \\
$(4,4\ ;\ 15,14,4,5,2,3,11,10,12,13,7,6,8,9,1,0)$ \\
$(4,4\ ;\ 15,14,6,7,11,10,2,3,12,13,5,4,8,9,1,0)$ \\
$(4,4\ ;\ 15,7,13,10,9,14,11,1,12,4,3,6,8,2,5,0)$ \\
$(4,4\ ;\ 15,5,13,7,9,1,11,3,14,4,12,6,10,2,8,0)$ \\
$(3,3,4\ ;\ 3,6,10,0,15,17,1,18,14,16,2,19,13,12,8,4,9,5,7,11)$ \\
$(3,3,4\ ;\ 3,6,11,0,8,19,1,15,4,17,14,2,16,18,10,7,12,9,13,5)$ \\
$(3,3,4\ ;\ 6,7,10,8,13,17,0,1,3,15,2,19,16,4,18,9,12,5,14,11)$ \\
$(3,3,3,3\ ;\ 3,4,6,0,1,9,2,15,17,5,13,18,19,10,23,7,22,8,11,12,21,20,16,14)$ \\
$(3,3,4\ ;\ 3,6,10,0,8,14,1,16,4,18,2,17,13,12,5,19,7,11,9,15)$ \\
$(3,3,4\ ;\ 3,6,10,0,8,14,1,18,4,16,2,19,13,12,5,17,9,15,7,11)$ \\
$(3,3,4\ ;\ 3,6,7,0,16,19,1,2,10,12,8,14,9,18,11,17,4,15,13,5)$ \\
$(3,3,4\ ;\ 6,10,19,8,13,17,0,12,3,15,1,16,7,4,18,9,11,5,14,2)$ \\
$(3,3,4\ ;\ 6,7,10,8,9,13,0,1,3,4,2,19,16,5,17,18,12,14,15,11)$ \\
$(3,3,4\ ;\ 3,6,10,0,8,18,1,14,4,16,2,19,13,12,7,17,9,15,5,11)$ \\
$(3,3,3,3\ ;\ 3,6,12,0,9,16,1,18,23,4,20,22,2,19,15,14,5,21,7,13,10,17,11,8)$ \\
$(3,3,3,3\ ;\ 3,4,6,0,1,9,2,15,17,5,14,13,19,11,10,7,23,8,22,12,21,20,18,16)$ \\
$(3,3,3,3\ ;\ 3,6,12,0,9,16,1,10,18,4,7,22,2,20,15,14,5,21,8,23,13,17,11,19)$ \\
$(3,3,3,3\ ;\ 3,6,12,0,9,16,1,18,11,4,23,8,2,19,15,14,5,21,7,13,22,17,20,10)$ \\
\hline\hline
\end{tabular}
\end{center}
\caption{Conjectured Mom-4's not identified by SnapPea.}
\end{figure}
Figure 2. The notation used can be interpreted as follows: the
numbers before the semi-colon describe the type of ideal polyhedra
used to construct the manifold. For example, the first entry in the
figure has the numbers ``$3,3,4$'' to the left of the semi-colon; each
``$3$'' indicates an ideal triangular dipyramid, while each ``$4$''
indicates an ideal square dipyramid (i.e. an ideal
octahedron). Each ideal dipyramid has two ``polar'' vertices and
either three or four ``equatorial vertices''. Number the faces of all
the polyhedra sequentially in such a way that the faces ``north'' of
each equator are numbered before the faces ``south'' of each
equator. For example, in the first entry the first triangular
dipyramid has faces 0, 1, and 2 next to one polar vertex, and faces
10, 11, and 12 next to the other polar vertex. The next triangular
dipyramid has faces 3, 4, and 5 as well as faces 13, 14, 15, and the
square dipyramid has faces 6 through 9 and 16 through 19. (This
somewhat unintuitive numbering scheme was chosen for convenience when
writing the computer software for this part of the paper.) Then the
numbers to the right of the semi-colon form a permutation which
describes how to glue together the faces of the ideal polyhedra.
For example, in the first entry the string of numbers which begins
with ``3, 6, 8, 0, \ldots'' imply that face 0 is glued to face 3, face
1 is glued to face 6, and so on. Since
we are requiring ``polar'' vertices to be identified solely with other
``polar'' vertices, no other information is needed to reconstruct the
polyhedral gluing.

One additional point of information:
all but eight of the manifolds in the list satisfy
\(|\partial M|=2\); seven satisfy \(|\partial M|=3\) and one satisfies
\(|\partial M|=4\). Thanks to the timely assistance of Morwen
Thistlethwaite, the authors were able to positively identify all eight
of these manifolds:

\begin{figure}
\begin{center}
\includegraphics{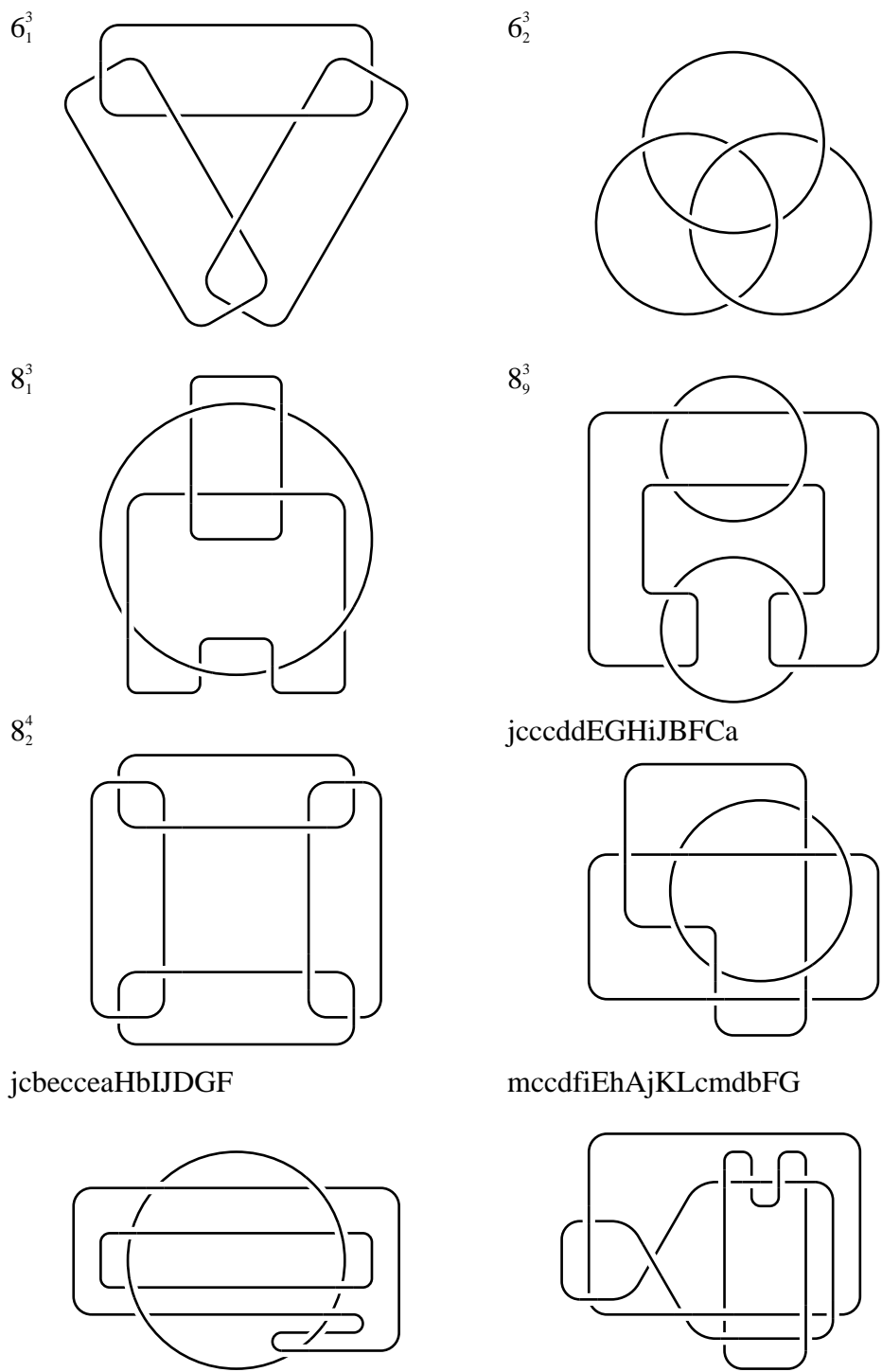}
\end{center}
\caption{Eight links whose complements are Mom-4's with 3 or more
cusps.}
\end{figure}

\begin{conjecture}
There are 8 hyperbolic manifolds $M$ such that \(M,T,\Delta\) is a
Mom-$n$ for some $2\le n\le 4$ and \(|\partial M|>2\). All eight
manifolds are complements of links in $S^3$: the links $6^3_1$,
$6^3_2$ (the Borromean rings), $8^3_1$, $8^3_9$, $8^4_2$, and the
links with Gauss codes jcccddEGHiJBFCa, jcbecceaHbIJDGF, and
mccdfiEhAjKLcmdbFG.
\end{conjecture}

At the time of writing we are still searching for an efficient way to
verify SnapPea's computations in the Mom-$4$ case; clearly, examining
over a million fundamental groups by hand is not a practical solution. Until
a better way is found, our enumeration results in the Mom-$4$ case should
properly be considered speculative.

\begin{section} {Hyperbolic Mom-$n$'s in  hyperbolic 3-manifolds are internal
Mom-$n$ structures for $n\le 4$}\end{section}

Let $R$ be a convolotube in the interior of a compact hyperbolic 
3-manifold $N$ and let $V$ be the cube
with knotted hole bounded by $R$.  By drilling out solid tori from 
$N-\cirV$, we can create a
manifold
$M $ which is non-elementarily embedded in $N$ and whose boundary 
contains a convolutube. We call such an embedding \emph{knotted}.
The goal
of this chapter is to show that if $n\le 4$, any  embedding
of a Mom-n manifold $(M,T)$ into a compact hyperbolic manifold 
$(N,T)$ is unknotted.

\smallskip
\begin{definition} Let $M$ be a compact 3-manifold and T a possibly 
empty union of components of
$\partial M$.  We say that
$(M,T)$ is
\emph{hereditarily unknotted}, if \emph{every} non-elementary
embedding into a  compact hyperbolic 3-manifold
$N$, taking $T$ to components of $\partial N$,  has the property that 
each component of  $\partial M$
  is either boundary parallel or bounds a solid torus.
\end{definition}

\begin{remark}  If $(M,T)$ is herditarily unknotted and
$M_1$ is obtained by filling a component of
$\partial M - T$, then $(M_1,T)$ is
hereditarily unknotted.\end{remark}

\begin{lemma}  If $(M,T)$ is a hereditarily unknotted Mom-$n$ manifold
non-elementarily embedded in the hyperbolic 3-manifold $N$ such that $T$ 
bounds a tubular neighborhood
of a geodesic, then $(M,T)$ is an internal Mom-$n$ structure.\end{lemma}

\begin{proof}  Let $V$ be the solid torus bounded by $T$.
By Lemma \ref{geodesic complement}, if $N_1=N-\cirV$
with cusp neighborhoods deleted, then $N_1$ is compact
hyperbolic.  Therefore
$(M,T)\subset (N_1,T)$ is a non-elementary embedding and
hence any  component of $\partial M-T$ either bounds a solid
torus or is boundary parallel in $N_1$.  Therefore similar
properties hold in $N$ and hence $(M,T)$ is an internal Mom-$n$
structure on
$N$.\end{proof}

\begin{remark} The condition that $T$ bounds a neighborhood of a 
geodesic is essential. \end{remark}

\begin{lemma}  \label{two}Let $M$ be a compact hyperbolic 3-manifold 
with $T$ a union of components of
$\partial M$.  If $\partial M-T$ is connected,  $(M,T)$ is hereditarily
unknotted.\end{lemma}

\begin{proof}  If under a non-elementary embedding $(M,T)\to (N,T)$,\ 
\ $\partial M-T$ was a
convolutube, then $M$ would be reducible.\end{proof}

The following result establishes criteria for showing that $(M,T)$ is 
hereditarily unknotted.

\begin{lemma}  \label{H criteria} Let $M$ be a compact  hyperbolic
3-manifold with $V_1$, \ldots, $V_n$ components of $\partial M$ and $T$ 
a nonempty union of some other
components.  If any of the following hold, there exists no 
non-elementary embedding $(M,T)\to
(N,T)$ such that $N$ is compact hyperbolic and $\{V_1,\ldots, V_n\}$ is exactly
the set of convolutubes of $\partial M\subset N$.

i) The manifold obtained by some filling of $M$ along $V_1$, \ldots,
$V_n$ is a 3-manifold without any hyperbolic part.  (That is, after
applying sphere and torus decompositions there are no hyperbolic
components.)

ii) After some filling of $M$ along $V_1$, \ldots, $V_n$, the surface $T$
is compressible.

iii) For every filling on a non-empty set of components of $\partial
M-T\cup V_1\cup \cdots\cup V_n$, either $ V_1\cup\cdots\cup V_n$ is 
incompressible or the
filled manifold has no hyperbolic part.
\end{lemma}

\begin{proof}   Suppose that $(M,T)$ embeds in $(N,T)$, where
among the components of $\partial M$, 
$V_1,\cdots, V_n$ are the set of convolutubes and $W_1,\cdots, W_m$
are the tubes.  Let $W_i^*$ denote the solid torus bounded by $W_i$ 
and $V_i^*$ denote the cube with
knotted hole bounded by $V_i$.  Let $B_1,\cdots,B_n$ be pairwise disjoint 3-balls 
in $N$ such that for each
$i$, $V_i\subset B_i$.

i) Let $\hat M$ be a manifold obtained by filling the $V_i$'s.  Let 
$\hat N$ be obtained by deleting the $V_i^*$'s and doing the 
corresponding fillings along the $V_i$'s.  Therefore $\hat N$ is 
obtained from $\hat M$ by Dehn filling and $\hat N$ is a connected 
sum of $N$ with $S^2\times S^1$'s and/or lens spaces and/or $S^3$'s. 
This implies that $\hat M$ has a hyperbolic part.

ii) If $T$ is compressible  in $\hat M$ it is compressible in $\hat 
N$ and hence in $N$, which is a contradiction.

iii) First observe that $V_i$ compresses in the manifold $M$' 
obtained by filling $M$  where each
$W_i$ is filled with $W_i^*$.  Topologically, $M^\prime$ is 
homeomorphic to $N$ with $n$ open
unknotted and  unlinked solid tori removed and so has a hyperbolic part.
\end{proof}

\begin{theorem}
If the Mom-$n$ manifolds for $n\le 4$ with three or more boundary
components are exactly those listed in Figure 3 (i.e. if Conjecture
\ref{list of mom-2's and mom-3's} is true), then any hyperbolic Mom-$n$
manifold $(M,T)$ with $n\le 4$ is hereditarily unknotted.
\end{theorem}

\begin{proof}  By Lemma \ref{two} it suffices to consider the case 
where $M$ is one of the eight Mom-4 manifolds with at least three
boundary components listed in Figure 3. If $M$ is any of the first six
manifolds and $T$ is any component of $\partial M$, then $(M,T)$ is
hereditarily unknotted by criterion (i) of Lemma \ref{H criteria}.
For manifolds 7 and 8, depending on which boundary component is used
for $T$, applications of (i) and (iii) imply that they are
hereditarily unknotted.
\end{proof}

\begin{section} {Examples of Mom-$n$ structures}\end{section}


In this section we give some representative examples of hyperbolic
manifolds $N$ which contain an internal Mom-2 or Mom-3 structure
\((M,T,\Delta)\).  Our goal in this section is to give the
reader an intuitive feel for how these particular cell complexes arise
inside hyperbolic manifolds. All of the manifolds in this section
involve manifolds $N$ with torus boundary, with the base torus of the
Mom-structure \((M,T,\Delta)\) being \(\partial N\).  To obtain Mom-$n$
structures on closed manifolds, note that if \(T=\partial N\) then a
Mom-$n$ structure \((M,T,\Delta)\) on $N$ passes to a Mom-$n$ structure on
any manifold obtained by filling \(\partial N\).


\smallskip
\begin{example}
The first example is the figure-$8$ knot complement. We construct a
Mom-2 \((M,T,\Delta)\) inside this manifold as follows. The torus $T$
is just the boundary of the manifold. The $1$-handles \(\lambda_1\)
and \(\lambda_2\) span the two tangles which make up the standard
diagram of this knot, as seen in
\begin{figure}[tb]
\begin{center}
\includegraphics{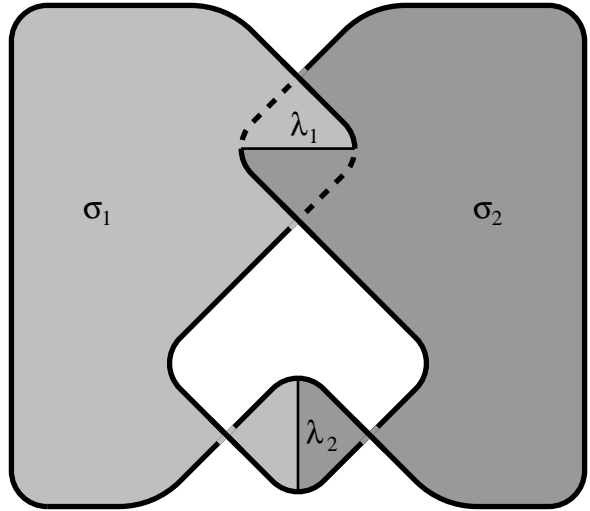}
\end{center}
\caption{The figure-8 knot complement equipped with a Mom-2.}
\end{figure}
figure 4. Finally the $2$-handles \(\sigma_1\) and \(\sigma_2\) are
symmetrically placed as shown in the diagram. Note that, as required,
each $2$-handle meets three $1$ handles counting
multiplicity. Specifically, each $2$-handle meets \(\lambda_1\) twice
and \(\lambda_2\) once. Also, one can see from the diagram that the
complement of $T\cup \{\lambda_i\}\cup \{\sigma_j\}$ consists of a
solid torus, and that the solid torus retracts onto a
homotopically non-trivial simple closed curve (which is a geodesic in
$N$).  Thus this is a valid hyperbolic Mom-2 structure on $N$.

Moreover, we can quickly determine the nature of the ideal
triangulation of $M$ described in the Section \ref{enumeration}. The
ends of \(\lambda_1\) are each dual to a four-sided pyramid in this
triangulation, and the two endpoints of \(\lambda_2\) are each dual to
a ``digonal pyramid'', so that each get eliminated. Thus the figure-8
knot complement possesses a Mom-2 structure \((M,T,\Delta)\) where $M$
is a two-cusped hyperbolic manifold which is in turn obtained by
gluing together the faces of an ideal octahedron. By the comments
after Theorem \ref{list of mom-2's and mom-3's}, $M$ must be the
complement of the Whitehead link.  And indeed, it is easy to verify
that if one drills out the core of the solid torus in the complement
of $M$ one obtains a manifold homeomorphic to the complement of the
Whitehead link.\end{example}


\begin{example}
Next we will let $N$ be the manifold known as m003 in the SnapPea
census. This manifold has first homology group \(\BZ + \BZ/5\), and
hence is not a knot complement; instead, we will present this manifold
as the union of two regular ideal hyperbolic tetrahedra; see
\begin{figure}[tb]
\begin{center}
\includegraphics{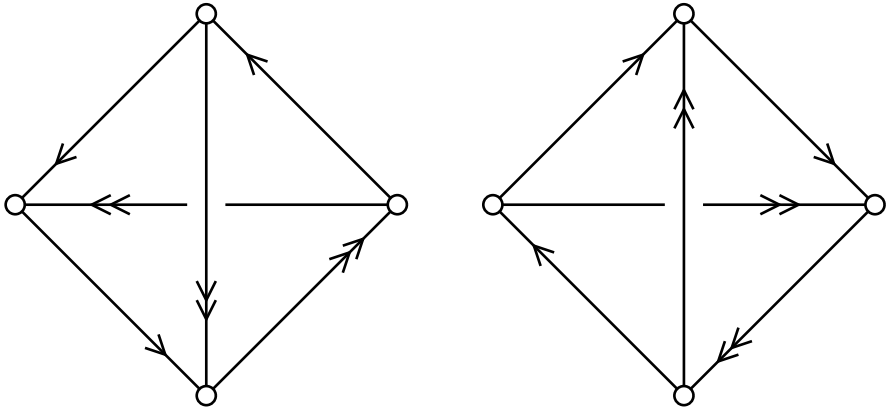}
\end{center}
\caption{The two ideal tetrahedra making up the manifold m003.}
\end{figure}
figure 5. Note that in the diagram each face is glued to the
corresponding face on the other tetrahedron, in such a way that the
edges match up into two equivalence classes as shown. To make $N$ a
compact manifold with torus boundary, assume the ideal tetrahedra are
truncated.  Now suppose that we construct $(M,T,\Delta)$ in this case
as follows.  For the $1$-handles, we use neighborhoods of the two
edges shown in the diagram, truncated by the torus \(T=\partial N\).
And for the $2$-handles, we use neighborhoods of the two truncated
triangles which are formed by gluing together the faces on the front
of each tetrahedron in the diagram. It is a simple exercise to confirm
that the complement of the resulting embedded manifold $M$ consists of
a solid torus, and that the solid torus retracts onto a simple closed
geodesic curve, and that therefore this manifold possesses a valid
hyperbolic Mom-2 structure. Each of the $1$-handles in this Mom-2
meets three of the $2$-handles, counting multiplicity; therefore we
can conclude that m003 contains a Mom-2 \((M,T,\Delta)\) where $M$ is
obtained by gluing together two ideal three-sided dipyramids. From
Theorem \ref{list of mom-2's and mom-3's} and the comments
following it we know this must be either m125 or m203. Further
investigation with SnapPea shows that it must in fact be m125.

It is instructive to get another view of this Mom-2 by constructing a
cusp diagram for this manifold. Specifically, consider the
triangulation induced on $T$ by the given ideal triangulation of
m003. The two ideal tetrahedra in m003 will appear as eight triangles,
the four ideal triangles will appear as twelve edges, and the two
edges will appear as four vertices. The resulting cusp diagram is
shown in
\begin{figure}[tb]
\begin{center}
\includegraphics{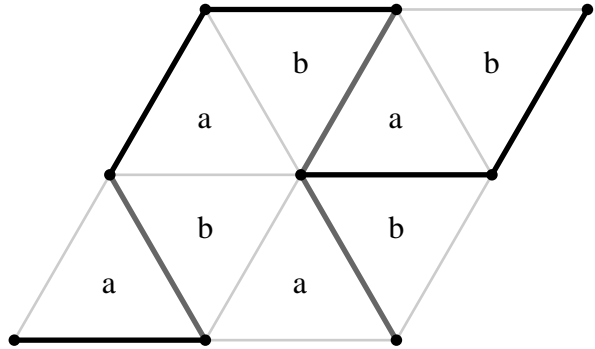}
\end{center}
\caption{The cusp diagram for m003, with the components of the Mom-2
highlighted.}
\end{figure}
figure 6; keep in mind this is a diagram of a torus, so the edges of
the parallelogram are identified with one another. (The labels inside
each triangle indicate which of the ideal simplices contributes that
triangle to the cusp diagram.)  The highlighted edges in the cusp
diagram are those that correspond to the $2$-handles of the handle
structure $\Delta$; in other words, they along with the four vertices
of the diagram comprise \(\Delta^1 \cap T\).\end{example}

\begin{example}
As another example in this vein, consider the manifold $N$=m017. This
manifold has first homology group \(\BZ+\BZ/7\), so again it is not a
knot complement in \(S^3\), but for brevity's sake we only present a
cusp diagram here. In
\begin{figure}[tb]
\begin{center}
\includegraphics{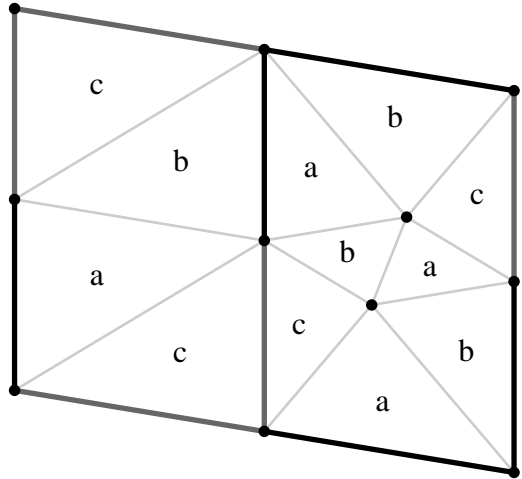}
\end{center}
\caption{The cusp diagram for m017, with the components of a Mom-2
highlighted.}
\end{figure}
figure 7, the corners of the three ideal hyperbolic tetrahedra which
make up m017 can be seen. And again, the highlighted edges in the
cusp diagram correspond to two faces of those tetrahedra which provide
the $2$-handles for an internal Mom-2 in this manifold. Note that we can
determine from the cusp diagram alone that the $1$-handles of this
Mom-2 meet four and two $2$-handles respectively, counting
multiplicity, and that therefore in the resulting Mom-2 structure
\((M,T,\Delta)\) the manifold $M$ is obtained by gluing together an
ideal octahedron. As before, this implies that $M$ must be
homeomorphic to the complement of the Whitehead link. Some further
work with SnapPea confirms this: m017 is obtained by (-7,2) Dehn
surgery on either component of the link.\end{example}

\begin{example}Finally, we include the motivating example for this paper.
\begin{figure}[tb]
\begin{center}
\includegraphics{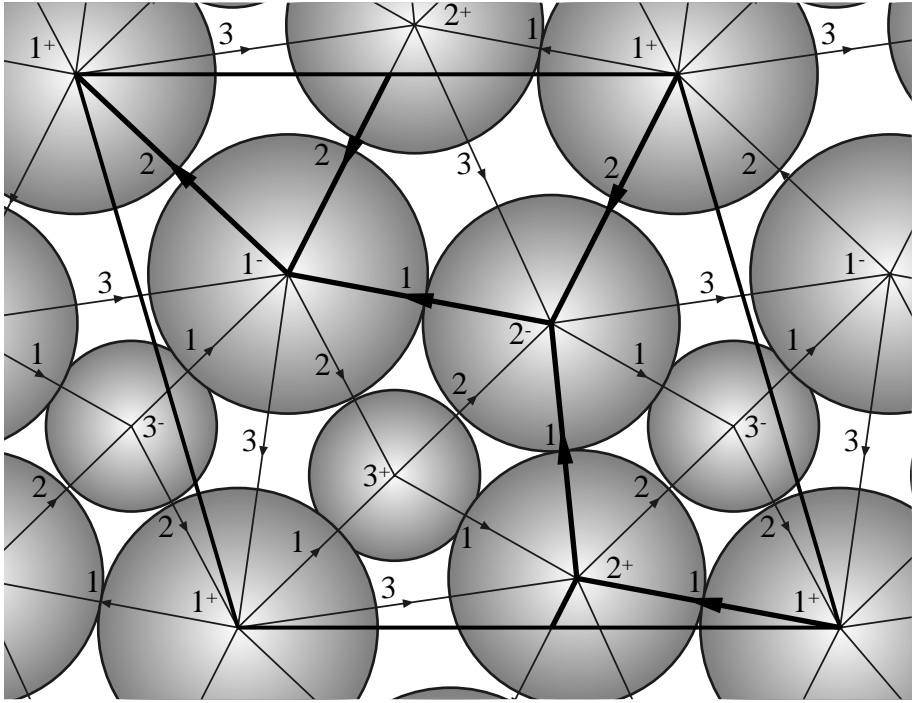}
\caption{SnapPea's cusp diagram for m011, with the components of a Mom-2
highlighted.}
\end{center}
\end{figure}
Figure 8 shows the maximal cusp diagram of the 1-cusped manifold m011
as provided by Weeks' SnapPea program. Unlike the previous cusp diagrams
in this section, it also shows all the horoballs at
hyperbolic distance at most 0.51 from the maximal horoball at
infinity. The paralellogram shows a fundamental domain for the
$\BZ\oplus\BZ$ action.  Note that the ideal triangulation presented in
this diagram is dual to the Ford decomposition of the manifold.
In particular the 1-simplices of the triangulation are geodesics
orthogonal to pairs of horoballs; these 1-simplices appear either as
edges in the figure joining the endpoints of the simplex in $\Sinfty$, or
as the vertical geodesics passing from the center of each horoball to the
horoball at infinity.

Let $H_\infty$ denote the horoball at infinity. There are six
horoballs in the diagram up to the $\BZ\oplus\BZ$-action, labelled
$1-$, $1+$, $2-$, $2+$, $3-$, and $3+$. This notation
means that if $\gamma\in \pi_1(m011)$ takes horoball $n\pm$ to
$H_\infty$, the horoball at infinity, then $H\infty$ is transformed to
one labelled $n\mp$.  The geodesic from $n-$ to $H_\infty$ is oriented
to point into $H_\infty$ and hence the geodesic from $n+$ to
$H_\infty$ is oriented out of $H_\infty$.  These orientations induce,
via the $\pi_1(m011)$-action, the indicated orientations on the edges
of the diagram.  We explain, by example, the meaning of the edge
labels.  The edge $3$ from $2+$ to $2-$ corresponds to a geodesic
$\sigma$ with the property that when $2+$ is transformed to
$H_\infty$, then $2-$ is transformed to $3+$ and $\sigma$ is
transformed to the vertical geodesic oriented from $H_\infty$ to $3+$.
(Had the edge been oriented oppositely, then $2-$ would have been
transformed to $3-$.)  SnapPea did not provide the orientation
information, however such information can be derived from the SnapPea
data.

By staring at this diagram we can see how m011 contains a Mom-$2$.  Let
$V_0$ be the maximal horotorus neighborhood of the cusp, slightly
shrunken.  By expanding $V_0$, the expanded $V_0$ touches the
(expanded) horoballs labeled $1$, thereby creating a 1-handle denoted
$E_1$.  Let $V_1$ denote this expanded $V_0$.  Further expansion
creates $V_2$ which is topologically $V_1$ together with another
1-handle $E_2$, this 1-handle occuring between horoball 2 and
$H_\infty$.  The edge labelled 1 between horoballs $2-$ and $2+$
corresponds to a valence three 2-handle which goes over $E_1$ once and
$E_2$ twice.  Similary the edge labelled $2$ between horoballs $1-$
and $1+$ gives rise to a valence three 2-handle going twice over $E_1$
and once over $E_2$.  The parallelogram of Figure 8 can also be viewed
as $\partial V_0$, with the centers of $1-,1+, 2-, 2+$ as the
attaching sites of the 1-handles and the thick black lines
corresponding to where the 2-handles cross over $\partial V_0$.
\end{example}

\vfill \pagebreak

\end{document}